\documentclass[a4paper]{article}

\usepackage{geometry}
\usepackage[utf8]{inputenc}
\usepackage[T1]{fontenc}
\usepackage[english]{babel}
\usepackage{lmodern}

\usepackage{indentfirst}
\usepackage{amsmath}
\usepackage{amssymb}
\usepackage{amsthm}
\usepackage{bm}
\usepackage{graphicx}
\usepackage[svgnames]{xcolor}
\usepackage[babel]{csquotes}
\usepackage{titlecaps}
\Addlcwords{a an the and or nor but as at by for in of on per to vs}

\usepackage{doi}
\newcommand{\email}[1]{\href{mailto:#1}{\texttt{#1}}}
\usepackage{hyperref}

\usepackage{thmtools}

\DeclareMathOperator{\e}{e}
\newcommand{\R}{\ensuremath{\mathbb{R}}}

\newcommand{\N}{\ensuremath{\mathcal{N}}}
\newcommand{\Z}{\ensuremath{\mathbb{Z}}}

\newcommand{\Y}{\mathbb{H}}
\renewcommand{\d}{\mathrm{d}}

\newcommand{\dy}{\d y}
\newcommand{\dz}{\d z}
\newcommand{\del}{\partial}
\renewcommand{\O}{\mathcal{O}}

\newcommand{\ra}{\rightarrow}

\newcommand\p{\partial}
\newcommand\mm{\mbox}
\newcommand\ve{\varepsilon}

\DeclareMathOperator{\tr}{tr}
\DeclareMathOperator{\diam}{diam}
\DeclareMathOperator{\vol}{Vol}
\DeclareMathOperator{\dvol}{\d\!\vol}

\DeclareMathOperator{\id}{Id}

\DeclareMathOperator{\sech}{sech}
\DeclareMathOperator{\length}{length}

\renewcommand{\epsilon}{\varepsilon}
\renewcommand{\Phi}{\varPhi}
\renewcommand{\Omega}{\varOmega} 

\newcommand{\diffp}[3][]{\frac{\partial^{#1}#2}{\partial {#3}^{#1}}}

\renewcommand{\le}{\leqslant}
\renewcommand{\ge}{\geqslant}
\newcommand{\delG}{\del_{\mathrm{G}}} 
\newcommand{\delM}{\del_{\mathrm{M}}} 

\DeclareMathOperator{\dist}{dist}
\newcommand{\eucl}{\ensuremath{\mathrm{Eucl}}}
\newcommand{\dd}{\textit{\dj\hspace{1pt}}}
\newcommand{\tdd}{\textit{\dh\hspace{1pt}}}
\newcommand{\red}{\mathcal{R}}

\newcommand{\thmref}[1]{\hyperref[#1]{\autoref*{#1} \enquote{\nameref*{#1}}}}
\newcommand{\thmnameref}[1]{\hyperref[#1]{\nameref*{#1} \ref*{#1}}}
\newcommand{\thmonameref}[2]{\hyperref[#1]{#2 \ref*{#1}}}

\newtheoremstyle{plaincaps}
    {}{}%
    {\itshape}{}%
    {\bfseries}{}%
    { }%
    {\thmname{#1}\thmnumber{ #2}\thmnote{ (\titlecap{#3})}}
\theoremstyle{plaincaps}
\newtheorem{theorem}{Theorem}[section]
\newtheorem{proposition}[theorem]{Proposition}
\newtheorem{defprop}[theorem]{Definition/Proposition}
\newtheorem{lemma}[theorem]{Lemma}
\newtheorem{corollary}[theorem]{Corollary}

\newtheoremstyle{definitioncaps}
    {}{}%
    {}{}%
    {\bfseries}{}%
    { }%
    {\thmname{#1}\thmnumber{ #2}\thmnote{ (\titlecap{#3})}}
\theoremstyle{definitioncaps}
\newtheorem{example}[theorem]{Example}
\newtheorem{examples}[theorem]{Examples}

\newtheorem{definition}[theorem]{Definition}

\theoremstyle{remark}

\hypersetup{%
	bookmarksopen,
	bookmarksnumbered,
	hidelinks,
	colorlinks,     
	linkcolor=DarkGreen,
	citecolor=Navy,
	urlcolor=Navy,
	pdftitle={Potential Theory on Gromov Hyperbolic Manifolds of Bounded Geometry},
	pdfauthor={Matthias Kemper, Joachim Lohkamp}
}

\begin{document}

\title{Potential Theory on Gromov Hyperbolic Manifolds of Bounded Geometry}
\author{Matthias Kemper, Joachim Lohkamp}
\date{}

\maketitle

\begin{abstract}
In the 1980s Alano Ancona developed a profound
potential theory on Gromov hyperbolic manifolds of bounded geometry. Since then, such hyperbolic spaces have become basic in geometry, topology and group theory. This also stimulated a broader interest into analytic aspects of these spaces. In this paper we make Ancona's original work, addressed to a rather advanced audience, approachable for a wider readership without firm background in potential theory. We streamline the arguments to make the results depend only on coarse geometro-analytic invariants and we derive useful details not found in the original sources. Finally, we explain remarkable relations of these results to recent developments in the potential theory and quasi-conformal geometry of domains in $\R^n$.
\end{abstract}

\tableofcontents

\section{Introduction}
In this paper we are interested in the elliptic analysis on \emph{non-compact} complete manifolds and its relation to the geometry at infinity. Non-compact spaces occur even if one is primarily interested in compact cases, for instance as limit spaces under blow-ups or by deforming  bounded domains. An example are Riemann uniformizations of planar domains which also turn the boundary of the domain into a boundary at infinity.

It is reasonable to focus on non-compact manifolds of \emph{bounded geometry}. In simple terms, bounded geometry excludes a variety of awkward analytic phenomena resulting from degenerations towards infinity. For instance,  it makes sure that basic analytic controls from local estimates, most prominently Harnack inequalities, hold uniformly all over the manifold. Moreover, a complete manifold of bounded geometry is  \emph{stochastically complete}. This means that solutions of stochastic differential equations associated to typical elliptic problems on the manifold exist for an infinitely long time.

However, having bounded geometry alone is not yet enough to gain a control as strong as in the classical case of harmonic or holomorphic functions on smoothly bounded planar domains. Finer details in the style of Cauchy's integral formula over the boundary of the given domain are reserved to the subclass of Gromov hyperbolic manifolds. For reasonable elliptic operators there is a contour integral for (positive) solutions formally taken over the so-called Martin boundary which, and this is a deep insight of Ancona's potential theory on such hyperbolic spaces, simply \emph{equals} the natural geometric boundary at infinity, namely the Gromov boundary.

\subsection{Goals of this paper}

Our main goals are to explain how bounded geometry paves the way to a controllable potential theory for suitable elliptic operators and how adding Gromov hyperbolicity refines these basic controls to a remarkable degree. We will see how even mild violations of Gromov hyperbolicity can cause a breakdown of these refined controls.

This is also a tour through a variety of beautiful concepts and fundamental results from geometric hyperbolicity notions over the potential theoretic balayage and the resolvent equation to (boundary) Harnack inequalities and their associated Martin theory.

The exposition is modelled on the original sources \cite{anc87} and \cite{anc90}, but we streamline and replace some of the arguments and derive a priori estimates only depending on coarse geometric and analytic data such as the hyperbolicity and coercivity constants of the considered spaces and operators. Our presentation emphasizes the central role of the boundary Harnack inequality on Gromov hyperbolic spaces and we derive some useful versions of this inequality not present in earlier literature.

Our second theme is a strategy, also initiated by Ancona  \cite{anc87,anc90}, to transfer the potential theory on hyperbolic spaces to Euclidean domains, provided they admit conformal deformations into such a hyperbolic space. The idea is to transform the operator in such a way that the notion of a solution remains unchanged under the deformation. The effect then is that the estimates and results derived on the hyperbolic space simply remain valid under this conformal deformation.

More recent work in quasi-conformal geometry due to Bonk, Heinonen and Koskela \cite{bhk01} has shown that so-called \emph{uniform domains} are precisely those domains admitting the type of  hyperbolizing  deformations needed here. Remarkably, the particular potential theory on uniform domains obtained this way also characterizes their uniformity condition. Aikawa  \cite{aik01,aik04} has developed an independent analytic argument showing that, under some weak regularity assumptions, uniform domains  are exactly those Euclidean domains which admit boundary Harnack inequalities for the Laplace operator. In turn, this underlines the sharpness of Ancona's results for the Gromov hyperbolic case.
\subsection{Organization of the paper}
\textbf{\hyperref[sec:bc]{Section~\ref*{sec:bc}}} \, We first recall the concepts of manifolds of bounded geometry and of Gromov hyperbolicity. We also introduce the ideal boundary of such space at infinity, the Gromov boundary. Then we discuss some analytic counterparts of these geometric concepts, uniform ellipticity and coercivity of operators, which are necessary for Ancona's theory.
Finally, we introduce the basic potential theoretic notion of balayage which we use to globalize local potential theoretic considerations.

\textbf{\hyperref[sec:loc]{Section~\ref*{sec:loc}}} \,  We use the bounded geometry and uniform ellipticity to derive local estimates for Green's functions of elliptic operators on uniformly sized balls and their relation to Harnack inequalities. The original sources \cite{anc87,anc90} employ advanced analytic tools due to Stampacchia \cite{sta65}. We reorganize the arguments and base them on more common and elementary techniques.

\textbf{\hyperref[sec:glob]{Section~\ref*{sec:glob}}} \,  Here we additionally use the coercivity properties to formulate resolvent equations for Green's functions considered as operators. We apply resolvent equations to
derive local and, by means of chains of Harnack inequalities, also global growth estimates for Green's functions.

\textbf{\hyperref[sec:hy]{Section~\ref*{sec:hy}}} \, Now we also use that our manifold is Gromov hyperbolic. We employ this hyperbolicity in two places, most notably in \thmref{prop:expg} to derive our main technical result, \thmref{thm:gphi}.
Secondly, the hyperbolicity readily implies the existence of an abundance of $\Phi$-chains which allows us to deduce the \thmnameref{cor:hbhi} and that the so-called Martin boundary is homeomorphic to the Gromov boundary. This gives us contour integrals for positive solutions of $Lw=0$ similar to the classical case of harmonic functions on smoothly bounded planar domains.

\textbf{\hyperref[sec:b]{Section~\ref*{sec:b}}} \,  Manifolds of bounded geometry that are Gromov hyperbolic form a rather special class of non-compact manifolds.
However, many other geometries admit canonical deformations into such hyperbolic manifolds. Basic examples are the so-called \emph{uniform} domains in $\R^n$. This remarkable concept covers smoothly bounded domains but also some classes  of fractals like snowflakes. We show how we can use the hyperbolic theory to gain control over the classical potential theory on such Euclidean domains from the insights of \autoref{sec:hy} .

\section{Basic Concepts}\label{sec:bc}

\subsection{Geometric Structures}\label{sec:gs}

We start with a discussion of those properties we assume throughout for our spaces and operators.
Our space $(M^n,g)$ is a connected, complete $C^\infty$-smooth Riemannian manifold of dimension $n \ge 2$. We assume the following two conditions. The first one, of having bounded geometry, is basic for a fine control of the potential theory. It is a homogeneity  condition saying that the space locally looks the same around any point, up to a uniformly controlled deviation:

\begin{description}
\item[Bounded Geometry] \emph{For each ball $B_{\sigma}(p) \subset M$ there is a smooth $\ell$-bi-Lipschitz chart $\phi_p$ to an open subset of $U_p$ of $\R^n$ with its Euclidean metric, where $\sigma={\sigma_M} >0$, $\ell=\ell_M \ge 1$ are fixed constants for the given manifold $M$. We call this a space of \emph{bounded geometry}, or more precisely, \emph{$(\sigma,\ell)$-bounded geometry}.}
\end{description}

On  spaces of bounded geometry one can link any two points $p,q \in M$ by a sequence of such balls transferring information inductively from one end to the other. Such configurations are known as Harnack chains due to their frequent use in the context of Harnack inequalities in elliptic analysis.

\begin{description}
\item[Harnack Chains] \emph{For some fixed $r \in (0,\sigma)$  we call a sequence of balls $B_r(x_1),\dots, B_r(x_k)$ with
\[\text{$x_1=p$, $x_k=q$ and $d(x_i,x_{i+1})<r/2$ for $i=1,\dots,k-1$}\]
a \textbf{Harnack chain} of length $k$.}
\end{description}
Thus any point $x_i$ is contained in the ball $B_{r/2}(x_{i+1})$. If we link $p$ and $q$ by a curve $\gamma_{p,q}$ with $\length(\gamma_{p,q})=d(p,q)$ we can cover $\gamma_{p,q}$ by balls forming a Harnack chain. From this we see that, for given $r$, we can choose a chain of length proportional  to the distance of the endpoints.

\begin{figure}[h]
\centering
\includegraphics[width=0.84\textwidth]{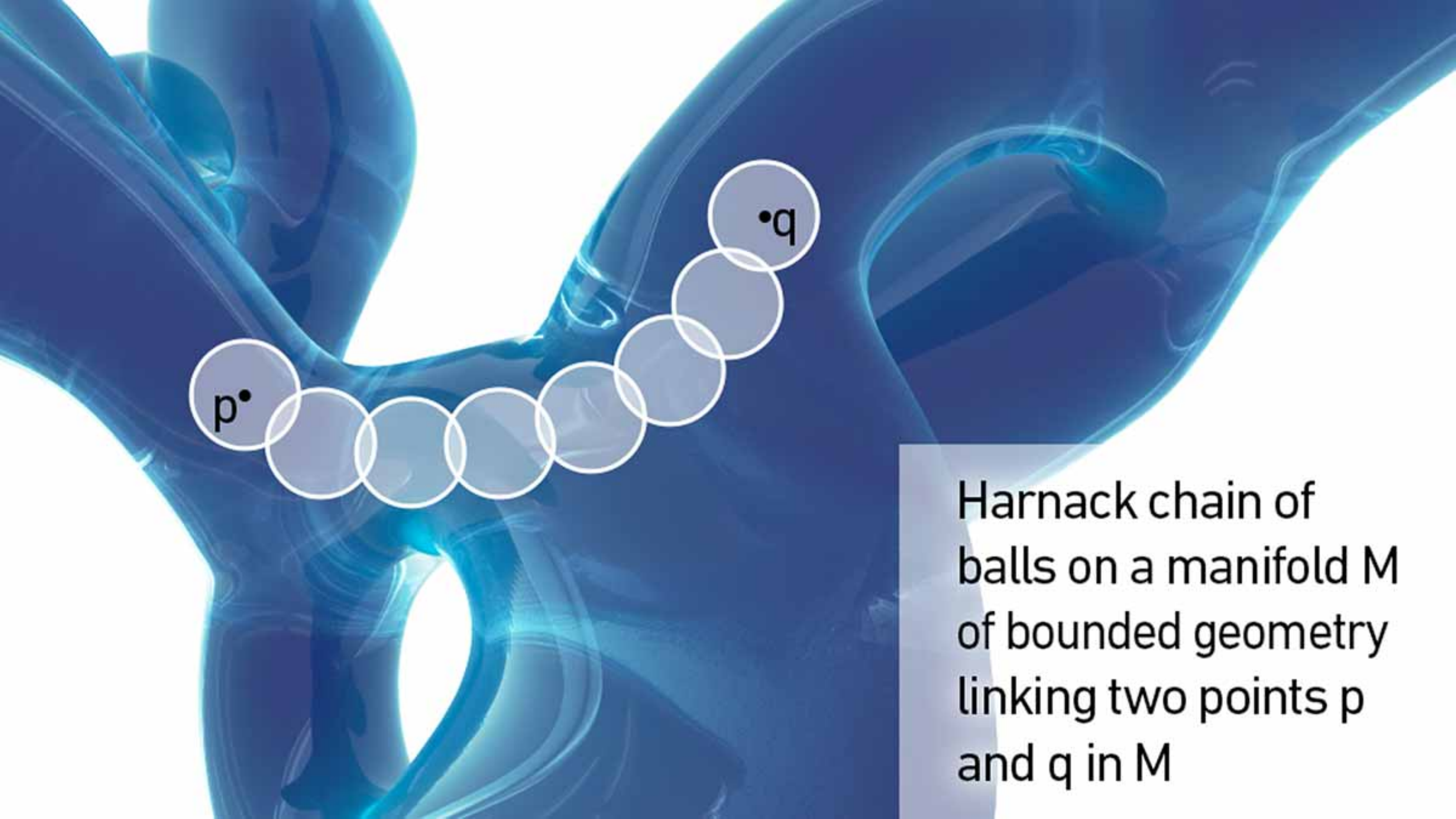}
\caption{Harnack Chains}
\end{figure}

The second condition, of Gromov hyperbolicity, that we impose on $M$ is, in some sense, complementary to the bounded geometry assumption. It has no local impact but strong consequences for the geometry near infinity.

\begin{description}
\item[Gromov Hyperbolicity]  \emph{A Riemannian manifold $(M^n,g)$, now regarded as a metric space, is Gromov hyperbolic or, quantitatively, \textbf{$\bm{\delta}$-hy\-per\-bo\-lic}, if all its geodesic triangles are \textbf{$\bm{\delta}$-thin}, for
some $\delta=\delta_M >0$. This means that each point on the edge of any geodesic triangle is within $\delta$-distance of one of the other two edges.}
\end{description}

The concept of Gromov hyperbolic spaces embraces a broad range of metric spaces, prototypes from either side are classical hyperbolic spaces $\mathbb{H}^n$ and abstract trees. While the former can be generalized to simply connected manifolds of sectional curvature less than a negative constant, $\delta$-hyperbolicity of a graph expresses that it is near to a tree on a large scale. Gromov hyperbolicity is designed to study the asymptotic behavior of spaces near infinity, cf.\ \cite[ch.~III.H]{bh99} and \cite{kb02} for detailed expositions.

Every compact manifold is Gromov hyperbolic since we can choose $\delta=\diam(M)$. This shows that this type of hyperbolicity has no local impact. However, near infinity Gromov hyperbolicity can be more demanding even than having constant negative sectional curvature. An instructive example are $\Z^2$-coverings of genus $\ge 2$ Riemann surfaces equipped with some hyperbolic metric. They have bounded geometry but they are not Gromov hyperbolic since we get arbitrarily large geodesic triangles uniformly close to Euclidean triangles in Hausdorff metric.

The motivating example for Ancona in \cite{anc87} were Cartan-Hadamard manifolds of sectional curvature pinched between two negative constants, which also ensures bounded geometry.
A larger and arguably more basic class of further examples for Gromov hyperbolic spaces of bounded geometry is given by uniform domains in Euclidean space with their quasi-hyperbolic metric as explained in \autoref{sec:b}.

\begin{description}
\item[Gromov boundary]\label{def:gprod}
\emph{For a complete Gromov hyperbolic space $X$ we define the set $\delG X$ of equivalence classes $[\gamma]$ of geodesic rays,  from a fixed basepoint $p \in X$, where two rays are equivalent if they have finite Hausdorff distance. $\delG X$ is called the \textbf{Gromov boundary} of $X$. ($\delG X$  does not depend on the choice of the basepoint $p$.)}
\end{description}

To define a topology on $\overline{X}^{\mathrm{G}} = X \cup \delG X$, we first recall that a \emph{generalized geodesic ray}  $\gamma: I \ra X$ is an isometric embedding of an interval $I \subset \R$ into $X$, where either $I = [0,\infty)$, then $\gamma$ is a proper geodesic ray, or $I = [0,R]$, for some $R \in (0,\infty)$. Then $\gamma$ is a geodesic arc. When we fix a basepoint $p \in X$ we can assign to any $x \in X$ a (not necessarily unique) generalized ray $\gamma_x$ with endpoint $\gamma_x(R) = x$. For the following discussion we extend the definition of such a ray to $I = [0,\infty]$ setting $\gamma_x(t) = x$ for $t \in [R,\infty]$.

Then the topology on $\overline{X}^{\mathrm{G}}$ can be characterized from the following notion of a converging sequence: $x_n \in \overline{X}^{\mathrm{G}}$ converges to $x \in \overline{X}^{\mathrm{G}}$ if there exist
generalized rays $\gamma_{x_n}$ with $\gamma_{x_n}(0) = p$ and $\gamma_{x_n}(\infty) = x_n$ subconverging (on compact sets) to a generalized ray $\gamma_x$ with $\gamma_x(0) = p$ and $\gamma_x(\infty) = x$. The canonical map $X \hookrightarrow \overline{X}^{\mathrm{G}}$ is a homeomorphism onto its image, $\delG X$ is closed and $\overline{X}^{\mathrm{G}}$  is compact, see \cite[III.H.3.7]{bh99}.
$\overline{X}^{\mathrm{G}}$ is called the \textbf{Gromov compactification} of $X$. It is a metrizable space.

To exploit Gromov hyperbolicity analytically Ancona introduced the following concept of $\Phi$-chains. As in the case of bounded geometry versus hyperbolicity, Harnack chains and $\Phi$-chains
act complementary. This is best seen in the  discussion preceding \thmref{prop:expg} where Harnack chains give basic estimates which, however, weaken the overall control, while $\Phi$-chains can be used to recover the apparently lost details.

\begin{description}
\item[$\Phi$-Chains]
\label{def:phi}
\emph{For a monotonically increasing function $\Phi:[0,\infty)\to(0,\infty)$ with $\Phi_0:=\Phi(0)>0$ and $\Phi(d)\overset{d\to\infty}{\longrightarrow}\infty$, a \textbf{$\Phi$-chain} on $X$ is a finite or infinite sequence $U_1\supset U_2\supset\cdots\supset U_m$ of open subsets of $X$ together with a sequence of \textbf{track points} $x_1,x_2,\dots,x_m$ such that
\begin{enumerate}
\item $\Phi_0\le d(x_i,x_{i+1})\le 3\Phi_0$,
\item $x_i \in \del U_i$,
\item $d(x,\del U_{i\pm1})\ge\Phi(d(x,x_i))$, for every $x\in \del U_i$
\end{enumerate}
for every $i$ where applicable.}\footnote{For notational convenience this is slightly different from Ancona's version in \cite[définitions V.5.1]{anc90}, but essentially the same.}
\end{description}

Note that a $\Phi$-chain traversed backwards, i.e., with sets $X\setminus U_m\supset X\setminus U_{m-1}\supset\cdots\supset X\setminus U_1$, is again a $\Phi$-chain with the same track points.
\begin{figure}[h]
  \centering
  \includegraphics[width=.84\textwidth]{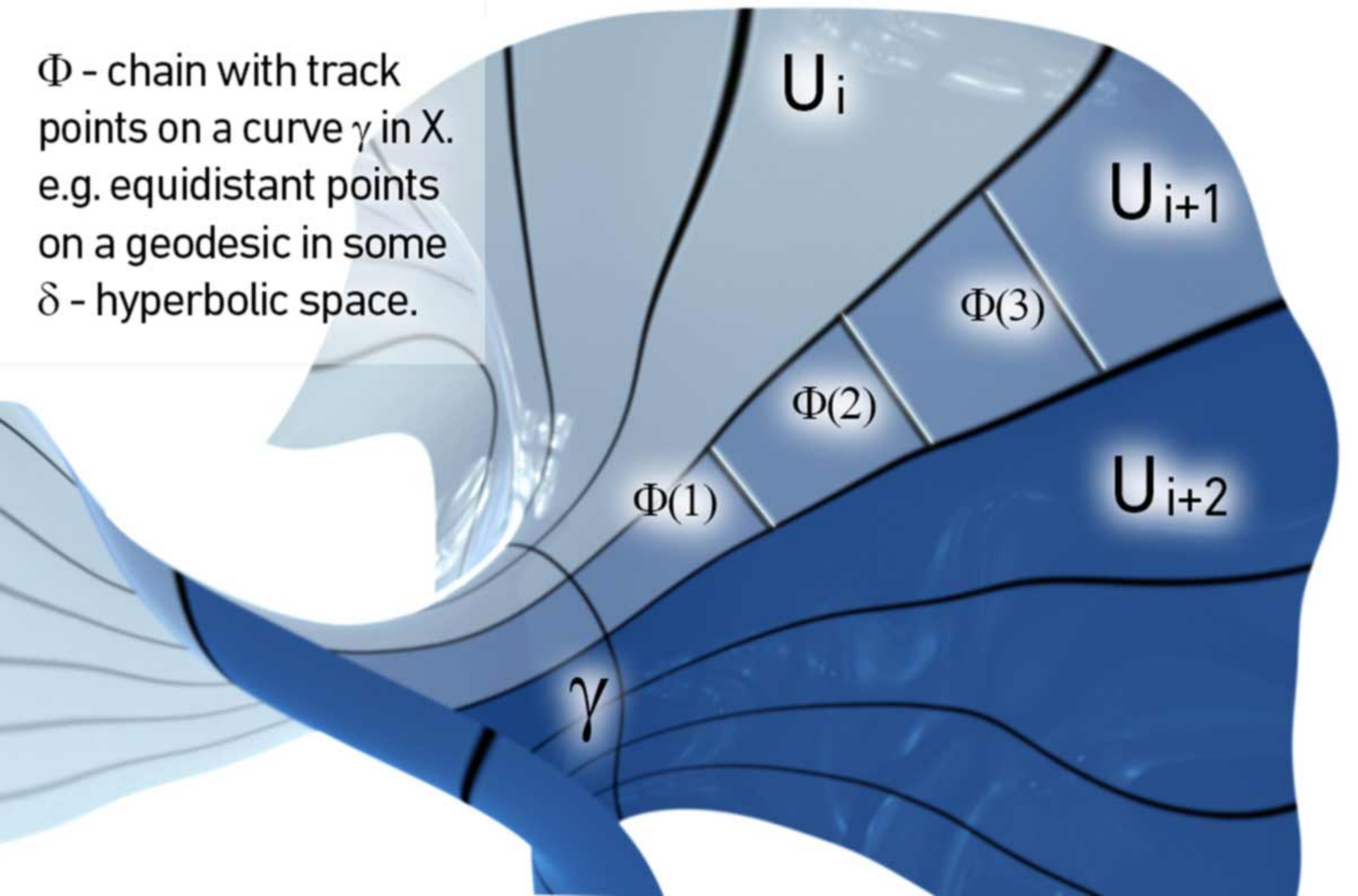}
  \caption{$\Phi$-Chains}
\end{figure}

The existence of infinite $\Phi$-chains can be considered as a partial hyperbolicity property of the underlying space. It is easy to see that neither the Euclidean nor asymptotically flat spaces admit any infinite $\Phi$-chains.

A first non-trivial example can be created as follows: with the coordinates $(x,y) \in \R \times \R^{n-1}$ we consider the metric $(1+|y|^2)^2 \cdot g_{\R} +g_\eucl$ on $\R \times \R^{n-1}$. Then the half-spaces $U_i:=\R^{\ge i} \times \R^{n-1}$ form a $\Phi$-chain with track points $x_i=(i,0)$ for $\Phi(t):=1+t^2$. But, as in the Euclidean case, the half-spaces  $U_i[k]:=\R^{k} \times \R^{\ge i} \times \R^{n-k-1}$, for $1 \le k \le n-1$, do not make up a $\Phi$-chain.

The hyperbolic $n$-space $\Y^n$ carries $\Phi$-chains in all directions and for the same $\Phi$.
Namely, when we think of $\Y^n$ as the upper half-space and $B_{1/2^i}(0)$ is the Euclidean ball of radius $1/2^i$, then $U_i=B_{1/2^i}(0)$, $i \ge 1$, form a $\Phi$-chain for $\Phi(t)=a + b \cdot t$ for suitable $a, b > 0$ independent even of $n$.  Due to the homogeneity of $\Y^n$ this also gives a $\Phi$-chain along any hyperbolic geodesic $\gamma$, i.e., with track points on $\gamma$.

This ubiquity of $\Phi$-chains which we see already from considering hyperbolic half-spaces relative to geodesics extends to arbitrary non-homogenous Gromov hyperbolic spaces.

\begin{description}
\item [$\Phi$-Chains on Hyperbolic Spaces]\cite[Section 8]{bhk01} \emph{For a general $\delta$-hyperbolic $X$ we choose a geodesic $\gamma: \R \ra  X$, and define \[ U_i:=\{x \in X\,|\, \dist\left(x,\gamma([22\delta i,+\infty))\right)<\dist\left(x,\gamma((-\infty,22\delta i])\right)\}\ .\] Then the $U_i$ form a $\Phi_\delta$-chain with track points $x_i = \gamma \cap \p U_i$ for $\Phi_\delta(t)=a_\delta + b_\delta \cdot t$ for  suitable $a_\delta, b_\delta  > 0$. Moreover, the $U_i$ and the $X \setminus U_i$ form \emph{neighborhood bases} in the topology of $\overline{X}^{\mathrm{G}}$ around the two ideal endpoints of the geodesic $\gamma$ in $\delG X$.}
\end{description}
The constants $a_\delta$ and $b_\delta$ can be inferred from the proof of \cite[Lemma 8.9]{bhk01}.

\subsection{Analytic Structures}
\label{sec:analytic}

Now we turn to the analytic side. On a complete Riemannian manifold $M$ with bounded geometry, we consider an operator $L$ with the following properties:
\begin{description}
\item[Adaptedness of $\bm{L}$] \emph{$L$ is a linear second order elliptic operator on $M$ such that relative to the charts $\phi_{p}$,
\[L(u) = -\sum_{i,j} a_{ij} \cdot \frac{\p^2 u}{\p x_i \p x_j} + \sum_i b_i \cdot \frac{\p u}{\p x_i} + c \cdot u
\]
for $(2,\beta)$-Hölder continuous $a_{ij}$, $(1,\beta)$-Hölder continuous $b_i$ and $\beta$-Hölder continuous $c$, $\beta \in (0,1]$, with
\[
k^{-1}\sum_i \xi_i^2 \le \sum_{i,j} a_{ij}  \xi_i \xi_j\le k \sum_i \xi_i^2
\]
for any $\xi\in\R^n$
and
 \[|a_{ij}|_{C^2(U_p)}, |b_i|_{C^1(U_p)}, |c|_{L^\infty(U_p)} \le k\]
for some $k=k_L \ge 1$ and any $p \in M$.}
\end{description}
Solutions of the equation $L \,w=0$ are called \textbf{\emph{L}-harmonic} functions. We also consider $L$-sub- and superharmonic functions satisfying, in the smooth case, $L \,w \le 0$ or $L \,w \ge 0$, respectively. More general notions of $L$-sub/superharmonicity will be introduced later on.\\

The adaptedness assumption is a generalized form of uniform ellipticity. It is a compromise between generality and readability, trimmed to ensure that one can apply the usual maximum principles, all weak solutions are classical $C^{2,\beta}$-regular solutions, the coefficients are bounded and that the adjoint operator $L^*$ shares (even quantitatively) the same properties because it still has Hölder continuous coefficients bounded by $3k$, cf.~\cite[Chapter 8]{gt98}. For more general conditions, see \cite{anc87}.

Towards the potential theory of such operators, we need the existence of a \textbf{Green's function} $G$ for $L$. This is a function $G(x,y)>0$ on $M \times M$ which is singular on the diagonal $\{(x,x)\:|\:x \in M\}$ and satisfies the equation $L\,G(\cdot,y)=\delta_y$, in a distributional sense, where $\delta_y$ is the Dirac delta function with basepoint $y$. Put differently, for any $y\in M$,  $G(\cdot,y) >0$ is a (singular) function so that for a given function $f$, the function $u(x)=\int_M G(x,y) \, f(y) \dvol(y)$ is a solution of $Lu=f$, where $\dvol$ is the volume measure on $M$. Then a Green's function associated to $L^*$ is given by $G^*(x,y)=G(y,x)$.

The existence of  a Green's function is not always granted. As a quite typical example for our further discussion we consider the Laplacian $-\Delta$ on the flat $\R^n$. The operator  $-\Delta- \lambda \id$ admits a Green's function if and only if $\lambda \le 0$. Moreover, we will see later on that there is a fine distinction between the cases $\lambda =0$ and $\lambda <0$. The following assumption we make for our operator $L$ means that we are focussing on the second case.

\begin{description}
\item[Weak Coercivity] \emph{For some $t>0$,  $L_t:=L-t\id$ admits
a  Green's function $G^t$.}
\end{description}

The weak coercivity condition is generic as soon as a Green's function exists. More precisely there is a \emph{generalized principal eigenvalue} $\tau=\tau_L$ such that for all $t<\tau$ the Green's function $G^t$ exists, while for $t>\tau$ there is no globally defined positive $L_t$-harmonic function. In the borderline case $t=\tau$, the Green's function $G^{\tau}$ might or might not exist, but in any case there is a globally defined positive $L_{\tau}$-harmonic function \cite[sec.~4.3]{pin95}. Then weak coercivity amounts to $\tau>0$.

Another way to express weak coercivity is to say that there is a positive supersolution $u$ of the equation $Lu = 0$ with the quantitative estimate $Lu \ge \ve \cdot u$, for some $\ve>0$.

Note that with $L$ the operator $L_t$ is again adapted and for $t<\tau$ even weakly coercive.

\subsection{Potentials and Balayage}
\label{sec:balayage}

The sum of a Green's function and a positive solution of $Lu=0$ is again a Green's function. Therefore one usually focusses on the uniquely determined \emph{minimal} Green's function which does not admit a sum decomposition into another Green's function and a positive solution.
More generally, any Green's function is an \emph{\textbf{$\bm{L}$-superharmonic}} function, that is a lower semi-continuous function with values in $(-\infty,\infty]$ that is larger than the Dirichlet solution on any ball with the same boundary conditions and finite on a dense set. On smooth parts, this is equivalent to $Lu\ge0$. An $L$-superharmonic function $p\ge 0$ \emph{without} a minorizing positive $L$-harmonic function $h$, that is with $p \ge h >0$, is called a \emph{\textbf{potential}}. Minimal Green's functions are not only examples for potentials, but also the basic building blocks:

\begin{theorem}[integral representation of potentials]\cite[22.]{her62}, \cite[Thm.\,6.18]{hm14}
Every potential $p$ on $M$ can be represented by a unique (positive) Radon measure $\mu_p$ as
\[p(x)=G(\mu_p)(x):=\int_M G(x,y)\,\d\mu_p(y)\ .\]
The support of $p$ (i.e.\ the complement of the largest open set where $p$ is $L$-harmonic) equals the support of $\mu_p$.
\end{theorem}

The fact that every positive $L$-superharmonic function is uniquely representable as the sum of a potential and a positive $L$-harmonic function is known as the \emph{Riesz representation theorem}. Hence, to get an integral representation for positive \emph{$L$-superharmonic} functions, only the $L$-harmonic part is left. The corresponding measures are supported on the \emph{Martin boundary} which is defined exactly for this purpose. Its identification in terms of more common geometric boundaries is the subject of subsections \ref{sec:mb} and \ref{sec:nu}.

Later we want to control $L$-superharmonic functions along $\Phi$-chains. Here we shift the part of the defining measure supported in $M \setminus  U_i$ onto $\del U_i$ without changing the function on $U_i$. This strategy is called sweeping or, due to its French origin (Poincar\'{e}, Cartan), \emph{\textbf{balayage}}.

Concretely, for an $L$-superharmonic function $u\ge0$  on $M$ and a subset $A\subset M$ we consider
\[\red_u^A:= \inf \{ v \ge0\:|\:v \mbox{ is $L$-superharmonic on } M \mbox{ with } v \ge u \mbox{ on } A \}\ .\]
This is called the \emph{reduit} (reduced).
It enjoys the following properties we will need later:
\begin{itemize}
\item $\red_u^A$ is $L$-harmonic outside of $\bar A$ and equal to $u$ in $A$.
\item The reduit is always $L$-superharmonic.
\item If $A$ is relatively compact, $\red_u^A$ is a potential.
\item $\red_{\lambda u}^{A}=\lambda \red_{u}^{A}$ for a constant $\lambda\ge0$, if $u$ is a positive function.
\item $\red_{u+v}^{A}=\red_{u}^{A}+\red_{v}^{A}$ for functions $u$, $v$.
\item $\red_u^{A\cup B}\le \red_u^A+\red_u^B$ for sets $A,B\subset M$.
\item Denoting the reduit with respect to the adjoint operator $L^*$ of $L$ by ${}^*\!\red$, we have $\red_{G(\cdot,y)}^A(x)={}^*\!\red_{G(x,\cdot)}^A(y)$ \cite[I.5.1, p.\,19]{anc90}.
\item If $G(\mu)$ is a potential, $\red_{G(\mu)}^A(x)=\int_M\red_{G(\cdot,y)}^A(x)\,\d\mu(y)$ for any $x\notin\bar{A}$ \cite[Théorème 22.4]{her62}. 
\end{itemize}

For general sets $A$, it may happen that $\red_u^A$ is not lower semi-continuous, but it always admits a canonical regularization $\widehat{\red}_u^A$, the \emph{balayée} (sweeped). For open sets or in general outside of $\bar{A}$ the two concepts coincide.
We can even recover the classical Perron solution $u$ of the Dirichlet problem on a ball $B$ with (say, continuous) boundary value $f$ as
\[ u(x)=\red_f^{\del B}(x)\ .\]

Also useful in this context are global variants of the maximum principle.
\begin{theorem}[global maximum principle]\cite[p.\,429]{her62}\label{gmp}
\begin{enumerate}
\item If $u$ is $L$-superharmonic on an open set $V\subset M$, $u\ge0$ on $\del V$, and there is a potential $p$ such that $u\ge-p$, then $u\ge0$ on $V$.
\item Let $p$ a potential, $L$-harmonic on an open set $V$ and locally upper bounded near every point of $\del V$. If $u\ge p$ on $\del V$ for some positive $L$-superharmonic function $u$, then $u\ge p$ in all of $V$.
\end{enumerate}
\begin{proof}
For (i), the function $\bar u$ defined as $\min(u,0)$ on $V$ and $0$ on $M\setminus V$ is $L$-superharmonic and $\ge-p$. Now the supremum of the family $\{\text{$L$-subharmonic functions $\le \bar u$}\}$ is $L$-harmonic,
$\ge -p$, and $\le0$, hence by the definition of potentials it is 0 which implies $u\ge0$.
\medskip

(ii) follows from (i) by considering the function $u-p$.
\end{proof}
\end{theorem}

\section{Local Maximum Principles and Harnack Inequalities}\label{sec:loc}
We derive estimates for positive solutions and Green's functions of $L_t$ on uniformly sized balls $B_R:=B_{\sigma/\ell}:=B_{\sigma/\ell}(0)\subset\R^n$ in the adapted charts. For notational convenience, in this chapter, all balls are measured in the Euclidean distance. The required assumptions on Locally defined, Adapted operators with a Green's function will be called

\begin{description}
\item[Assumption (LAG)]\label{ass:lac} We say that a differential operator $L$ on $B_R(0)\subset\R^n$ fulfills assumption (LAG) for $R>0$, $k>0$ and $n\ge2$  
    if
    \begin{itemize}
    \item $L$ is \emph{$k$-adapted} on $B_R(0)$ and
    \item there is a Green's function for $L$ on $B_R(0)$.
    \end{itemize}
\end{description}
\def\reflac{\hyperref[ass:lac]{(LAG)}}
Later all results can be transferred to the manifold using the bi-Lipschitz charts with only a small loss in constants.
Weak coercivity will not yet be used explicitly, but the results and estimates apply to all operators $L_t$ with $0\le t\le\tau$, where the constants now depend on $\tau$ as well (via $k$), but not on $t$.

\subsection{Maximum Principles}
The following version of the weak maximum principle does not require positivity of the \enquote{potential} term $c$ in $L$, the existence of a Green's function (or equivalently, a positive $L$-superharmonic function) is sufficient, but the result is weaker. It can easily be obtained from the version for $c\ge0$ using the $h$-transform, see \cite{pin95}, where the sign convention for $L$ is different from ours.
\begin{theorem}[local minimum-zero principle]\cite[Theorem~3.2.2, p.\,81]{pin95}\label{lmzp}
If $L$ satisfies assumption \reflac, $u$ is continuous on $B_R$, $L$-superharmonic on a domain $D\Subset B_R$ and $u|_{\del D}\ge 0$, then $u\ge 0$ in all of $D$.
\end{theorem}

With help of the uniform bounds on the coefficients, this can be upgraded to quantitative bounds in terms of the boundary values.
\begin{lemma}[local almost-maximum principle]\label{lamp}
If $L$ satisfies \reflac, there are constants $0<r_\mathrm{am}<R$ and $m>0$ depending only on $k$ and $n$, such that any $L$-subharmonic function $u$ on $B_r$, $r<r_\mathrm{am}$, with $u|_{\del B_r}\le m$ satisfies $u\le 1$ on $B_r$.
\begin{proof}
In polar coordinates, $L$ applied to a radial function $f(r)$ can be written as
\[ Lf=-\alpha f''(r)+\frac{\alpha-\tr a}{r}f'(r)+bf(r)+cf(r) \]
with $k^{-1}\le\alpha:=\frac{a_{ij}(x)x^ix^j}{r^2}\le k$, $nk^{-1}\le\tr a\le nk$ and $|b|$, $|c|$ bounded by $k$.

Consider the function $f(r)=\sech(\beta r)$ with some $\beta>0$. We want to have the following properties:
\begin{itemize}
\item $Lf\ge 0$ for $r<r_\mathrm{am}$,
\item $f\le1$, and
\item $f|_{\del B_r}\ge\sech(\beta r_\mathrm{am})=:m$ for $r<r_\mathrm{am}$.
\end{itemize}
Only the first property is not obvious and we need to tune the parameters $r_\mathrm{am}$ and $\beta$ to achieve it.
We can calculate
\[ Lf=\left[\alpha \beta^2\left(1-\sinh^2(\beta r)\right)\sech^2(\beta r)+(\tr a-\alpha)\beta^2\frac{\tanh(\beta r)}{\beta r}-b\beta\tanh(\beta r)+c\right]\sech(\beta r)\ .\]
Here the first \emph{two} terms behave as $\beta^2$ for large $\beta$ and $\beta r<1$. But in leading order in $\beta r$ the $\alpha$ terms cancel, while the positive $\tr a$ term survives. Using the uniform bounds on $\alpha$ and $\tr a$ in terms of $k$ and $n$, we first fix a small $\beta r_\mathrm{am}$ so that the $\beta^2$ terms are positively bounded below and then make $\beta$ large to kill the remaining terms and get $Lf\ge0$ for $r<r_\mathrm{am}$ where $r_\mathrm{am}$ and $\beta$ depend only on $k$ and $n$.

Now we can apply the \thmnameref{lmzp}
to $f-u$, because $L(f-u)\ge 0$ and $(f-u)|_{\del B_r}\ge 0$. This yields $u\le f\le1$ on $B_r$.
\end{proof}
\end{lemma}

With completely different methods we can see that Green's functions behave for small distances like the Euclidean Green's function
\[
    G_\eucl(x,y)=\begin{cases}
                         \frac{1}{2\pi}\log|x-y|\ , & \text{if } n=2 \\
                         \frac{1}{(n-2)\vol S^{n-1}}\frac{1}{|x-y|^{n-2}}\ , & \text{if } n\ge3\ .
                       \end{cases}
\]
\begin{lemma}[local bound for some Green's function]\label{lbsgf}
For $L$ satisfying \reflac\ there are constants $r_\mathrm{gb}>0$, $\tilde q\ge1$ depending only on $k$ and $n$ and some (not necessarily minimal) Green's function $\tilde G(x,y)$ such that
\[
	\tilde q^{-1}\,G_\textup{eucl}(x,y)\le\tilde G(x,y)\le \tilde q\,G_\textup{eucl}(x,y)\quad\text{for }x,y\in B_{r_\mathrm{gb}}\ .
\]
\begin{proof}
We sketch the argument for $a_{ij}=\delta_{ij}$ and $n\ge3$. The general case is merely notationally more involved and can be found in \cite[pp.\,61--63]{mir70}.

As distributions, we have $L_xG_\eucl(x,y)=\delta(x-y)+R(x,y)$, where $R(x,y)$ has a singularity at $x=y$ of order $\O(r^{-(n-1)})$. This means that we can choose $r_\mathrm{gb}>0$ so that $\int_{B_{r_\mathrm{gb}}}|R(x,y)|\dy<1/2$ for $x\in B_{r_\mathrm{gb}}$ because the integral scales as $\O(r_\mathrm{gb})$ and all involved constants depend only on the bounds on the coefficients $k$. This in turn implies the finiteness of all iterates $R^n$ where the product of kernels is defined as $R^2(x,y)=\int_{B_{r_\mathrm{gb}}}R(x,z)R(z,y)\dz$. The bound ensures also the summability of
\[
    \tilde G:=G_\eucl\sum_{i=0}^\infty(-R)^n\ ,
\]
again as products of kernels.
Application of $L$ yields
\[
    L\tilde G=(\delta+R)\sum_{i=0}^\infty(-R)^n=\sum_{i=0}^\infty(-R)^n-\sum_{i=1}^\infty(-R)^n=\delta\ ,
\]
so that $\tilde G$ is indeed a Green's function. $\tilde G-G_\eucl$ has a singularity at $x=y$ of quantitatively lower order than $G_\eucl$, hence the explicit bounds.
\end{proof}
\end{lemma}

\begin{corollary}[local bound for the minimal Green's function]\label{lbmgf}
If $L$ fulfills \reflac, there are constants $q>0$ and $0<r_\mathrm{mgbi}<r_\mathrm{mgbo}:=\min(r_\mathrm{am},r_\mathrm{gb})$ depending only on $k$ and $n$ such that the \emph{minimal} Green's function $g(x,y)$ of $B_{r_\mathrm{mgbo}}$ (i.e., vanishing on the boundary $\del B_{r_\mathrm{mgbo}}$) satisfies
\[
    q^{-1}\,G_\textup{eucl}(x,y)\le g(x,y)\le q\,G_\textup{eucl}(x,y)\quad\text{for }x,y\in B_{r_\mathrm{mgbi}}\ .
\]
\begin{proof}
For fixed $y\in B_{r_\mathrm{mgbo}/2}$, $g$ can be represented as $g(x,y)=\tilde G(x,y)-u(x)$ where $u$ is the solution of the Dirichlet problem $Lu=0$ in $B_{r_\mathrm{mgbo}}$ and $u(x)=\tilde G(x,y)$ on $\del B_{r_\mathrm{mgbo}}$. By the \thmnameref{lbsgf} we know
\[
    u(x) \le\frac{\tilde q}{|x-y|^{n-2}}\le\frac{\tilde q}{(r_\mathrm{mgbo}/2)^{n-2}}\qquad\text{for }x\in\del B_{r_\mathrm{mgbo}}
\]
(and similarly for $n=2$). Thus we can apply the \thmnameref{lamp} to get $u\le\tilde q/(m(r_\mathrm{mgbo}/2)^{n-2})$, whence
\[
    g(x,y)\ge\frac{\tilde q^{-1}}{|x-y|^{n-2}}-\frac{\tilde q}{m (r_\mathrm{mgbo}/2)^{n-2}}\quad\text{for }x\in B_{r_\mathrm{mgbo}}\ .
\]
From this equation, $q\ge\tilde q$ and $r_\mathrm{mgbi}$ can be determined to obtain the lower bound.
For the upper bound, note that $\tilde G$ is positive and thus $g\le\tilde G$ by the \thmnameref{lmzp}.
\end{proof}
\end{corollary}

\subsection{Harnack Inequalities}
This standard result can be found e.g.\ in \cite[Theorem~8.20]{gt98}, but we include a proof here because it is an easy consequence of the previous results.
\begin{theorem}[Harnack inequalities]
\label{harnackineq}
For $L$ satisfying \reflac, there is an $H(k, n)\ge 1$ such that if $u>0$ is $L$-harmonic on $B_r(x_0)\subset B_R$, for some $0<r\le R$, then
\[ H^{-1}u(x_0)\le u(x) \le H u(x_0)\quad\text{for any } x \in B_{r/2}(x_0)\ .\]
\begin{proof}
We start with the case $r\le r_\mathrm{mgbi}$.

The potential $\red_u^{B_r(x_0)}$ (on the total space $B_{r_\mathrm{mgbo}}$) can be represented as $\red_u^{B_r(x_0)}=g(\mu_u)$ for some positive measure $\mu_u$ with support in $\del B_r(x_0)$. On $B_r(x_0)$, $\red_u^{B_r(x_0)}$ coincides with $u$ and we have $u=g(\mu_u)$.
Now for $y\in B_{r/2}(x_0)$, we can apply both estimates from \thmref{lbmgf} to get
\begin{multline*}
	u(y)=\int_{\del B_r(x_0)}g(y,z)\,\d\mu_u(z)\le \int_{\del B_r(x_0)} \frac{q}{(r/2)^{n-2}}\,\d\mu_u(z)\\\le 2^{n-2}q^2\int_{\del B_r(x_0)}g(x_0,z)\,\d\mu_u(z)=2^{n-2}q^2\,u(x_0)\ .
\end{multline*}
The other inequality is analogous.

For the case of larger $r$, we can employ the Harnack chains introduced in \autoref{sec:gs} with radius $r_\mathrm{mgbi}$. With every ball the estimate collects the constant for that case, but the total number of balls needed is bounded above by $R/r_\mathrm{mgbi}$, so we can just use the resulting power of $2^{n-2}q^2$ as $H$ in every case.
\end{proof}
\end{theorem}
The Harnack chain tactic just demonstrated will be used in the following to get estimates for even larger distances. The result are exponential bounds (with fixed constants) for the growth of $L$-harmonic functions---that is not especially satisfying, but a good starting point for improvement.

This Harnack inequality always needs a little more space around the balls where it holds, but for instance on bounded domains, it would be useful to have a similar result near points on the boundary. A blueprint is the following classical result. 
\begin{theorem}[boundary Harnack inequality on a disc]\label{thm:bhid}
There exist constants $A$, $C>1$ such that for any point $\xi\in \p B_1(0)\subset\R^2$ and $0<R<1$ the following is true: For any two harmonic functions $u$, $v>0$ (with respect to the Laplacian) on $B_{A \cdot R}(\xi) \cap B_1(0)$ that vanish along $B_{A \cdot R}(\xi) \cap \p B_1(0)$,
\[
\frac{u(x)}{v(x)} \le C   \frac{u(y)}{v(y)} \mm{ for all } x,\, y \in B_R(\xi) \cap B_1(0)\ .
\]
\end{theorem}
In the non-boundary version, the appearance of another function $v$ is obscured by the fact that the constant function 1 is (nearly) harmonic. Furthermore, the restriction to functions vanishing on the relevant part of the boundary is unavoidable, yet in combination with balayage techniques there are powerful applications as we will see after the proof of a much more general boundary Harnack inequality in \autoref{sec:bhp}.

\section{Global Results from Resolvent Equations and Boun\-ded Geometry}\label{sec:glob}
Now we use that our manifold $M$ has Bounded geometry and $L$ is Adapted and weakly Coercive, more precisely:
\begin{description}
\item[Assumption (BAC)]\label{ass:bac} We say that the pair $(M,L)$ fulfills assumption (BAC) for $\sigma>0$, $\ell\ge1$, $k>0$, $n\ge2$ and $\tau>0$ if $L$ is a differential operator on a connected complete noncompact Riemannian manifold $M^n$ such that
    \begin{itemize}
    \item $M$ is of \emph{$(\sigma,\ell)$-bounded geometry},
    \item $L$ is \emph{$k$-adapted} in the bounded geometry charts, and
    \item $L$ is \emph{weakly coercive} with generalized principal eigenvalue $\tau=:2\theta$.
    \end{itemize}
\end{description}
\def\refbac{\hyperref[ass:bac]{(BAC)}}

This gives us global growth estimates for the minimal Green's function $G$ of $L$. The basic idea is to combine the local estimates we derived in the last chapter with the resolvent equation for Green's functions viewed as operators.

While the boundedness of geometry and adaptedness of $L$ allow to carry over all results from \autoref{sec:loc} with constants now depending on $k$, $n$, $\sigma$, $\ell$ and $\tau$, as explained in the beginning of that section, the weak coercivity assumption is the most important new ingredient in this section.
Note that we do not yet make use of Gromov hyperbolicity.

\subsection{Resolvent Equation}

For a closed  operator $L$ defined on a dense set of a Banach space the set $\rho(L)$ of all $\lambda$, so that  the \emph{resolvent} $R_\lambda=(L-\lambda\id)^{-1}$ exists and is continuous, is called the \emph{resolvent set}. $\rho(L)$ known to be an open set. For $\lambda, \mu \in \rho(L)$ we have the resolvent equation \[R_\lambda \circ R_\mu= (\lambda  - \mu)^{-1} \cdot (R_\lambda - R_\mu).\]

One of the main applications of this identity in the context of elliptic operators is the comparison of solutions of $R_\lambda w=0$ with those of $R_\mu w=0$. We will use this idea in the special case of the minimal Green's functions $G^t$ corresponding to the operators $L^t=L-t\id$ on $M$, viewed as resolvents. We include the simple proof of the resolvent equation in our context and note some of its consequences.

\begin{lemma}[resolvent equation]
\label{resolvent} Assume $(M,L)$ satisfies assumption \refbac. For any $0\le t<\tau$, the minimal Green's function $G^t$, viewed as an operator on the space of positive Radon measures $\mu$ with $G^t\mu<\infty$, satisfies
\[G^t=G+tG\circ G^t\ .\]
Since all involved operators are positive, this yields the inequalities
\begin{align}
  G & \le G^t \quad\text{and}\label{eq:resineq1}\\
  G\circ G^t & \le \frac{1}{t}G^t\ .\label{eq:resineq2}
\end{align}
In particular, these results hold in the cases of characteristic functions of bounded measurable sets or Dirac measures.
\begin{proof}
Consider an increasing sequence of relatively compact, smoothly bounded open sets $(U_i)$ with $\bigcup U_i=M$. On each $U_i$, the corresponding Dirichlet Green's function $G_i$ (i.e., $G_i(\cdot,y)|_{\del U_i}\equiv0$) satisfies the resolvent equation in the form
\[ G_i(x,y)=G_i^t(x,y)-t\int_{U_i}G_i(x,z)G_i^t(z,y)\dz \]
because the right-hand side fulfills the properties of a Green's function (application of $L$ yields $\delta_y$) and has the correct boundary behavior. To see that the integral is finite for $x\neq y$, notice that any Green's function is integrable near its pole by the same arguments as in the proof of \thmref{lbmgf}.

For $i\to\infty$, the $G_i$ and $G_i^t$ converge to $G$ and $G^t$, respectively, uniformly on compact sets away from the pole because they are increasing and bounded above by the minimal Green's function. The same is true if the first argument is fixed because then one has the Green's function for the adjoint operator. It is easy to see that the equation survives in the limit $i\to\infty$.

By integration and Fubini's theorem one obtains the resolvent equation for arbitrary Radon measures $\mu$, as soon as $G^t\mu<\infty$ (and thereby $G\mu<\infty$).
\end{proof}
\end{lemma}

Recall from \autoref{sec:analytic} that unlike $L_t$, for $t < \tau$, the operator $L_{\tau}$ is no longer weakly coercive. In turn, we observe from  \eqref{eq:resineq2} that using resolvents loses its strength when $t$ approaches $0$.  This suggests to work with  $L_\theta$, for
\[\theta:=\tau/2.\]

\subsection{Behavior of Green's functions}
We combine the resolvent equation with \thmnameref{harnackineq} to derive growth estimates of  minimal Green's functions. 
\begin{proposition}[bound for the Green's function]
\label{prop:bgreen}
Given $(M,L)$ satisfying \refbac, there is a constant $c_1(\sigma,\ell, k, n,\tau)\ge 1$ such that for the minimal Green's function of $L$,
\begin{align*}
c_1^{-1}\le\ &G(x,y) &&\text{if }d(x,y)\le\sigma\text{, and}\\
&G(x,y)\le c_1 &&\text{if }d(x,y)=\sigma\ .
\end{align*}
\begin{proof}
The lower bound is directly obtained from the \thmnameref{lbmgf} because we have $G(\cdot,y)\ge g(\cdot,y)$ for any Green's function $g$ on a smaller domain and iterated application of the \thmnameref{harnackineq} (a number depending on $\sigma$, $\ell$ and $r_\mathrm{mgbi}$).

For the upper bound, consider the ball $B:=B_{\sigma/3}(x)$. The function $G^\theta(\chi_B)=\int_{B}G^\theta(\cdot,\zeta)\d\zeta$ is \emph{bounded} on $B$, where $\chi_B$ denotes the characteristic function of $B$, because we can apply the lower bound $1\le c_1 \cdot G(\zeta,z')$ to some $z'$ with $d(x,z')=2\sigma/3$ and arrive at
\[\int_{B}G^\theta(z,\zeta)\d\zeta \le c_1\int_{M}G^\theta(z,\zeta)G(\zeta,z')\d\zeta \le  \frac{c_1}{\theta} G^\theta(z,z')\]
for all $z\in B$. In the second inequality we used \eqref{eq:resineq2} for the Dirac function in $z'$. By continuity of $G^\theta(\cdot,z')$ on the compact $\overline{B}$ this is uniformly bounded.

Thus we can integrate the \thmnameref{resolvent} applied to $\chi_B$ over $B$ and get
\[\infty>\int_M\chi_BG^\theta(\chi_B) = \int_M\chi_B G(\chi_B)+\theta\int_M\chi_BG(G^\theta(\chi_B))\ . \]
When $G^*$ denotes the minimal Green's function for the adjoint operator $L^*$ this means
\[0 <\int_M\chi_BG(\chi_B) = \int_MG^\theta(\chi_B)\left[\chi_B-\theta G^*(\chi_B)\right]\ .\]
Hence, there must be an $x'\in B$ such that $\left(\chi_B-\theta G^*(\chi_B)\right)(x')>0$ and therefore
\[ \int_BG(\zeta,x')d\zeta=G^*(\chi_B)(x')< 1/\theta \]
which in turn shows the existence of an $x''\in B$ with $d(x',x'')>\sigma/12$ and
\[G(x'',x')\le\theta^{-1}\left(\vol(B\setminus B_{\sigma/12}(x'))\right)^{-1}\ .\]
Now the bounded geometry assumption assures that the volume can be bounded below by a constant $V_{\sigma,\ell,n}$ depending only on $\sigma$, the Lipschitz constant $\ell$ and the dimension $n$.

Fom the bounded geometry constraint, we observe that a Harnack chain of at most $\lceil24+\pi\ell\rceil$  balls $B_r (p_i)$, $i=1,...,\lceil24+\pi\ell\rceil$, of radius $r=\sigma/12$ with $d(p_i, p_{i+1})<r/2$  and $d(p_i, x'')>r$ suffices to link $x'$ with either $x$ or $y$.  We get similar chains starting from $x''$ to $y$ or $x$, respectively. Applying the \thmnameref{harnackineq} on each of these balls and multiplying the Harnack constants step by step to get $H^i$ on $B_r (p_i)$ we get the assertion for $c_1=H^{2 \cdot \lceil24+\pi\ell\rceil} \cdot \theta^{-1}V_{\sigma,\ell}^{-1}$.
\end{proof}
\end{proposition}
Note that we could now easily get more explicit growth estimates for $G$ by comparison with the \thmnameref{lbmgf}, but for the following these coarse estimates are sufficient.

\subsection{Relative Maximum Principles}
We start with an elementary comparison result which already visualizes the central effect of weak coercivity: with the same boundary conditions, compared to $L_\theta$-harmonic functions, $L$-harmonic functions \emph{sag}, imagine a rope\footnote{This analogy is quite precise, a suspended rope forms a hyperbolic cosine, a solution of the shifted one-dimensional Laplace equation.} or a rubber blanket where you apply less and less tension. To see this, we start with a paraboloid that fits in between.

\begin{lemma}
Under the local assumptions \reflac, there is a constant $m(R=\sigma/\ell,k,n)>0$ such that we can find a smooth function $f$ defined on $B_{R/2}(0)$ with
\[ f(x)=-m\ ,\quad f|_{\del B_{R/4}(0)}\ge 0\quad\text{and}\quad Lf>-1\ .\]
\begin{proof}
In polar coordinates, $L$ applied to a radial function $f(r)$ can be written as
\[ Lf(r)=-af''(r)+\frac{b}{r}f'(r)+cf(r) \]
with functions $a$, $b$, $c$ bounded by a constant depending only on $k$, $n$ and $R$.

With the ansatz $f(r)=s(r^2-\rho)$ we have $Lf(r)=-2as+2bs+cs(r^2-\rho)$ and first choose $\rho>0$ sufficiently small to assure $f\ge 0$ outside of $B_{R/4}(0)$ and then make $s>0$ small to achieve $Lf\ge -1$. This only uses the bounds on coefficients such that $f(0)=-s\rho=:-m<0$ depends only on $k$, $n$ and $R$.
\end{proof}
\end{lemma}

\begin{proposition}[relative maximum principle, local version]
\label{qmploc}
Under assumption \reflac\ for $L$ and $L_\theta=L-\theta\id$, assume further we have two functions $u$ and $\bar u>0$ on $B_R=B_R(0)$, $u$ $L$-subharmonic ($Lu\le0$) and $\bar u$ $L_\theta$-harmonic on $B_R$ with $\bar u|_{\del B_{R/4}}\ge u|_{\del B_{R/4}}$.
Then there is a constant  $\tilde\eta=\tilde\eta(R, k,  n,\tau) \in (0,1)$ such that
\[ u(0)\le \tilde\eta\,\bar u(0)\ .\]
\begin{proof}
Consider the function $h(z):=\bar u(z)+\theta f(z)\inf_{w\in B_{R/4}}\bar u(w)-u(z)$ where $f$ is the function from the previous lemma on $B_{R/4}$.
On $B_{R/4}$, we have $L\bar u=L_\theta\bar u+\theta\bar u\ge \theta\bar u$, $L f>-1$ and therefore
\[Lh(z)>\theta\left(\bar u(z)-\inf_{w\in B_{R/4}}\bar u(w)\right) \ge 0\ .\]
The boundary condition together with $f|_{\del B_{R/4}}\ge 0$ yields $h|_{\del B_{R/4}}\ge 0$ and we can apply the \thmnameref{lmzp}
for $L$ to see $h\ge 0$ in the interior of $B_{R/4}$ and especially $h(0)\ge0$, with $f(0)=-m$ we have
\[ u(0)=\bar u(0)-mt\inf_{w\in B_{R/4}}\bar u(w)\le(1-m \theta H^{-1})\bar u(0)\ . \]
The harmonicity on the full ball $B_R$ was used only in the last step to apply the \thmnameref{harnackineq}.
\end{proof}
\end{proposition}

Applied globally, this describes the relative growth of  $L$-harmonic versus $L_\theta$-harmonic functions.

\begin{proposition}[relative maximum principle, global version]
\label{qmpglob}
If $(M,L)$ satisfies \refbac, there is a constant $\eta=\eta(\sigma,\ell,k,n,\tau)\in(0,1)$ such that the following holds:\\
Assume we have two functions $u$ and $\bar u$ defined on $B_{r+3}(x)$ for some $x\in M$ and $r\ge\sigma$, $u$ $L$-subharmonic and $\bar u$ $L_\theta$-harmonic on $B_{r+3}(x)$ and $\bar u|_{\del B_r(x)}\ge u|_{\del B_r(x)}$. Then
\[ u(x)\le \eta^r \bar u(x)\ .\]
\begin{proof}
For integer multiples $r$ of $\sigma/4\ell$, this is proven by inductively applying \hyperref[qmploc]{the local version} in adapted charts along a chain of intersecting balls of length proportional to $r/(\sigma/\ell)$. On each of them we may apply the same Harnack inequality as described in the proof of \ref{prop:bgreen} above, so that we can choose $\eta$ as a function of $\sigma$, $\ell$, $H$ and $\tilde\eta$.
\end{proof}
\end{proposition}

For Green's functions, we get the following variants. We do not use them in the following arguments but they are worth being mentioned since they give some non-trivial constraints on the Green's functions from
our standard assumption \refbac.

\begin{corollary}[exponentially stronger decay]
\label{prop:expstrdec} Under assumptions \refbac,
there are constants $A(\sigma,\ell, k,  n,\tau)>0$ and $\alpha_1(\sigma,\ell, k,  n,\tau)>0$ such that
\[ G(x,y)\le A\e^{-\alpha_1 d(x,y)}G^\theta(x,y)\quad\forall x,y\in M\ . \]
\end{corollary}
From this we get the following growth estimate for Green's functions using the resolvent equation.
\begin{proposition}[exponential decay]
\label{prop:expdec} Under assumptions \refbac,
for suitable constants $B(\sigma,\ell, k,  n,\tau)>0$ and $\alpha_2(\sigma,\ell, k,  n,\tau)>0$ we have
\[ G(x,y)G(y,x)\le B \e^{-\alpha_2 d(x,y)}\quad\text{for }d(x,y)> 2\sigma\ .\]
\begin{proof}
Let $x'\in M$ such that $d(x,x')=\sigma$. Then employing the \thmnameref{resolvent} we have
\begin{multline*}
	G(x,y)G(y,x)\le G(x,y)\,H^3\,G(y,x')
\le \frac{H^5}{\vol(B_{\sigma/2}(y))}\int_{B_{\sigma/2}(y)} G(x,z)G(z,x')\d z\\
    \le c\int_M G(x,z)G(z,x')\d z
	\overset{\eqref{eq:resineq1},\eqref{eq:resineq2}}{\le} \frac{c}{\theta}G^\theta(x,x')\le\frac{cc_1^\theta}{\theta}
\end{multline*}
where we used the \thmnameref{harnackineq} for $L^*$ and the \thmnameref{prop:bgreen} applied to $L_\theta$ which itself satisfies the assumptions, but it may lead to a weaker constant $c_1^\theta$.

For the very same reason, we can do all of the above with $G^\theta$ instead of $G$ to get the uniform boundedness of $G^\theta(x,y)G^\theta(y,x)$, again with slightly worse constants. Combined with \thmref{prop:expstrdec} we have proved the assertion.
\end{proof}
\end{proposition}

In the case of a self-adjoint operator $L=L^*$ the Green's function is symmetric and the proposition says that it decays exponentially with the distance. This result does not use Gromov hyperbolicity and holds also e.g.\ in Euclidean space. The difference from the familiar Euclidean Laplacian's Green's function is owed to weak coercivity.

\section{Hyperbolicity and Boundary Harnack Inequalities}\label{sec:hy}

Now we additionally invest the hyperbolicity of the underlying space. The property we employ is that in Gromov hyperbolic spaces any two of their points can be connected by well-controlled $\Phi$-chains, cf.~\autoref{sec:gs}. The results we prove are more general. They hold with respect to any individual $\Phi$-chain even if the space carries essentially only this one $\Phi$-chain as in  the example $(\R \times \R^{n-1},(1+|y|^2)^2 \cdot g_{\R} +g_\eucl)$.

We exploit this to derive the main results for the potential theory on these spaces, the particular behavior of Green's function along $\Phi$-chains, building on the results for bounded geometries we derived in the last few sections. From this we infer the boundary Harnack inequalities which largely dominate the potential theory of our elliptic operators.

Our general assumptions \refbac{} on the manifold $M$ and the elliptic operator $L$ remain the same as in the previous section: $M$ is complete with bounded geometry and $L$ is adapted and weakly coercive. The additional assumptions, that is, the presence of a $\Phi$-chain (depending on the function $\Phi$) or even of an underlying hyperbolic geometry (with constant $\delta$ and coming with a universal function $\Phi=\Phi_\delta$), are stated directly in the results.

\subsection{Global Behavior: \texorpdfstring{$\Phi$}{Phi}-Chains}\label{subsection:gb}

The following result describes the key feature of $\Phi$-chains in the potential theory of our elliptic operators.

\begin{theorem}[Green's functions along $\Phi$-chains]
\label{thm:gphi}
Under assumptions \refbac, there is a suitable constant $c(\sigma,\ell, k,  n,\tau,\Phi) >1$ such that for any $\Phi$-chain with track points $x_1,\dots,x_m$ we have for the minimal Green's functions
\[c^{-1}G(x_m,x_j)\,G(x_j,x_1)\le G(x_m,x_1)\le c\,G(x_m,x_j)\,G(x_j,x_1)\ ,\quad j=2,\dots,m-1\ .\]
\end{theorem}

At the heart of the argument we employ the pairing of two at first sight entirely unrelated geometric and analytic properties: the existence of $\Phi$-chains on $M$  and the weak coercivity of $L$.
The idea is that, on the one hand, $\Phi$-chains  allow to find balls of arbitrary large radii in $U_{i-1} \setminus U_{i+1}$ centered in $\p U_i$ within a uniformly upper bounded distance to the track points. On the other hand, the \thmonameref{qmpglob}{relative maximum principle} shows that on these balls we can improve estimates we have diminished from a previous application of a
Harnack inequality. This makes the following result the main step in the proof of the Theorem.

\begin{proposition}[growth recovery along $\Phi$-chains]
\label{prop:expg}
With assumptions \refbac, for any given $\Phi$-chain with track points $x_1,\dots,x_m$ we have for the minimal Green's functions
\begin{equation}
\label{eq:step1}
G(z,x_1)\le c\,G^{\theta}(z,x_j)\,G(x_{j+1},x_1)\quad\text{for }z\in\del U_{j+1},
\end{equation}
for some constant $c(\sigma,\ell, k,  n,\tau,\Phi) >0$ independent of the length $j$ of the $\Phi$-chain.

\begin{figure}[h]
  \centering
  \includegraphics[width=.84\textwidth]{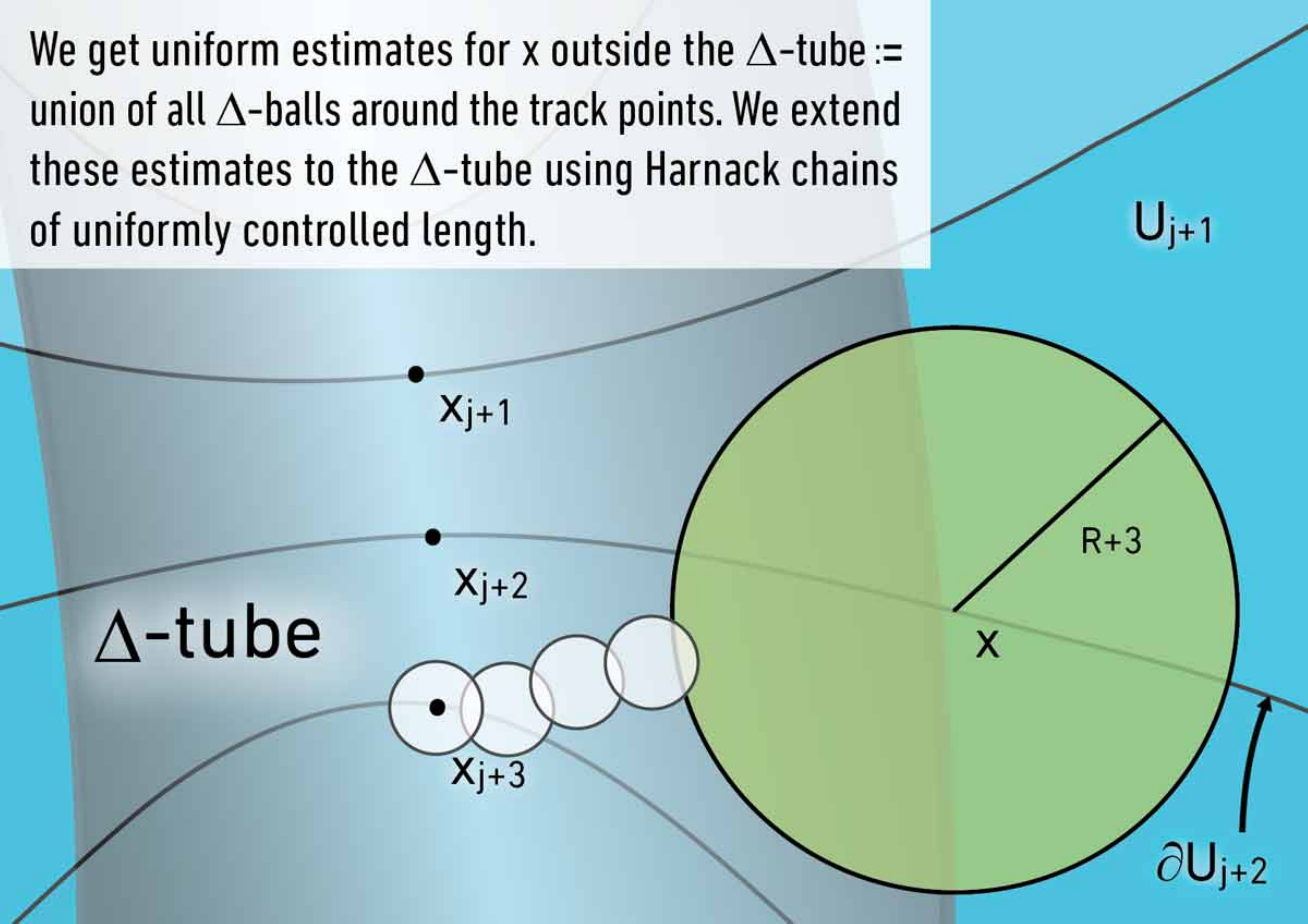}
  \caption{Growth Recovery Mechanism}
\end{figure}

\begin{proof} The argument is by induction over the length $j$. For $j=1$, the result follows from $G^\theta\ge G$, inequality  \eqref{eq:resineq1} and the \thmonameref{prop:bgreen}{lower bound for the Green's function}, $c_1\,G(x_2,x_1)\ge1$, from where we can take the first guess for the constant $c=c_1$ and note that $c_1$ depends only on $\sigma$, $\ell$, $k$, $\tau$ and $\Phi$.
For the induction step we first assume we have proved the weaker assertion that there is a constant $c_j$ so that \eqref{eq:step1} holds for any $\Phi$-chain of length $j$.  Then we can apply the Harnack inequalities for $L$ and $L_\theta^*$ to get a constant $c'(\sigma,\ell,k,\tau, \Phi)\ge 1$, independent of $j$, such that
\begin{equation}
\label{eq:step1a}
G(z,x_1)\le c'c_j\,G^{\theta}(z,x_{j+1})\,G(x_{j+2},x_1)\quad\text{for }z\in\del U_{j+1}\ .
\end{equation}
By the \thmnameref{gmp} this inequality extends to $z\in \bar U_{j+1}$. Now we invest the weak coercivity of $L$ and the properties of the $\Phi$-chains to improve this inequality.\medskip

Towards this end, we first apply the \thmonameref{qmpglob}{relative maximum principle} to the ($L$-superharmonic) function $G^{\theta}(\cdot,x_{j+1})$ and its greatest $L$-harmonic minorant $u$ on some ball $B_R(x)$ which we can represent as $u=\red_{G^{\theta}(\cdot,x_{j+1})}^{\del B_R(x)}$, the reduit always taken w.r.t.\ $L$. For $R=\ln(1/c')/\ln \eta$ and $B_{R+3}(x) \subset U_{j+1}$, the relative maximum principle yields
\[ u(x)\le \frac1{c'}G^{\theta}(x,x_{j+1})\ .\]

In turn, the definition of a $\Phi$-chain shows that there is some $\Delta(\Phi,R)>0$ such that $B_{R+3}(x)\subset U_{j+1}$, for $x\in\del U_{j+2}$, as soon as $d(x,x_{j+2})\ge\Delta$. Then, we have from  \eqref{eq:step1a}
\begin{align*}
G(x,x_1)&=\red_{G(\cdot,x_1)}^{\del B_R(x)}(x)\le c'c_j \red_{G^\theta(\cdot,x_{j+1})}^{\del B_R(x)}(x)\,G(x_{j+2},x_1)\\
&=c' c_j u(x)\, G(x_{j+2},x_1) \le  c_j G^{\theta}(x,x_{j+1})\, G(x_{j+2},x_1).
\end{align*}

From this we also get estimates for $x\in\del U_{j+2}$ with $d(x,x_{j+2})<\Delta$. This uses the \thmonameref{prop:bgreen}{bounds for the Green's function} and, as in the proof of \ref{prop:bgreen}, Harnack inequalities along chains of balls where length and radius are controlled by $\Delta$ (and the usual parameters for $M$ and $L$). We  get a constant $c''(\Delta, H)\ge1$ such that
\[ G(x,x_1)\le c'' c_j \,G^{\theta}(x,x_{j+1})\,G(x_{j+2},x_1)\quad\text{for $x\in\del U_{j+2}$ with $d(x,x_{j+2})<\Delta$}\ .\]
Now we compare this with  \eqref{eq:step1a}. The critical improvement we made is that for $d(x,x_{j+2})\ge\Delta$, independently of $j$, we get  \eqref{eq:step1a}
for $c_{j+1}=c_1$ and thus for  $c_{j+1}=c'' c_1$ in general. That is, \eqref{eq:step1} holds for $c=c'' c_1$.
\end{proof}
\end{proposition}

\begin{proof}[Proof of {\thmref{thm:gphi}}]
The first inequality is rather easy:
For $x\in\del B_\sigma(x_j)$ we have
\[ G(x,x_j)G(x_j,x_1)\le c_1 G(x_j,x_1)\le c_1 H\,G(x,x_1) \]
by \thmref{prop:bgreen} and the \thmnameref{harnackineq}. Since the left hand side is a potential and the right hand side is $L$-superharmonic, this inequality extends to $M\setminus B_\sigma(x_j)$ and in particular to $x_m$ by the \thmnameref{gmp}.

For the second inequality we use repeatedly the \thmref{prop:expg} and the \thmnameref{resolvent}:
\begin{align}
G(x_m,x_1)&=\red_{G(\cdot,x_1)}^{\del U_j}(x_m) &&|\ x_1\notin U_j\nonumber\\
 &\overset{\eqref{eq:step1}}{\le}c\,\red_{G^\theta(\cdot,x_j)}^{\del U_j}(x_m)\, G(x_{j+1},x_1)\nonumber\\
 &=c\,\red_{G(\cdot,x_j)+\theta G(G^\theta(\cdot,x_j))}^{\del U_j}(x_m)\,G(x_{j+1},x_1)  &&|\text{ res.eq.}\nonumber\\
 &\le c\left(G(x_m,x_j)+\theta \int_M \red_{G(\cdot,z)}^{\del U_j}(x_m)\,G^\theta(z,x_j)\,\dvol(z) \right)G(x_{j+1},x_1)\label{eq:long}
\end{align}
At this point we can again employ the first step \eqref{eq:step1}, but now for the reversed $\Phi$-chain $x_m,\dots,x_1$ with $M\setminus U_m,\dots,M\setminus U_1$ and for the adjoint operator $L^*$, namely
\[ G(x_m,z)\le c\,G^\theta(x_{j+2},z)\,G(x_m,x_{j+1})\quad\text{for }z\in M\setminus U_{j+1}\ .\]
This holds on all of $M\setminus U_{j+1}$ by the \thmnameref{gmp}. Since $x_j\in M\setminus U_{j+1}$, this can be directly applied to $G(x_m,x_j)$.
For the second summand in \eqref{eq:long}, we have $\red_{G(\cdot,z)^{\del U_j}}^{\del U_j}(x_m)={}^*\!\red_{G(x_m,\cdot)}(z)\le G(x_m,z)$ for $z\in \del U_j\subset M\setminus U_{j+1}$ (denoting the reduit with respect to $L^*$ by ${}^*\!\red$, cf.\ \autoref{sec:balayage}), but the upper bound ${}^*\!\red_{G(x_m,\cdot)}^{\del U_j}(z)\le c\,G^\theta(x_{j+2},z)\,G(x_m,x_{j+1})$ is valid for \emph{all} $z\in M$ by the definition of the reduit since the right hand side is positive and $L^*$-superharmonic in $z$.
Thus,
{\small \begin{multline*}
G(x_m,x_1) \le c^2\,G(x_m,x_{j+1})\,G(x_{j+1},x_1)\left(G^\theta(x_{j+2},x_j)+\theta\int_M G^\theta(x_{j+2},z)\,G^\theta(z,x_j)\,\dvol(z)\right)\ .
\end{multline*}}The large bracket is universally bounded by \thmref{prop:bgreen}, the \thmnameref{harnackineq}, and the inequalities \eqref{eq:resineq1} and \eqref{eq:resineq2} following from the \thmnameref{resolvent} for $t=\frac32\theta$.
\end{proof}

Now we assume that $M$ is a $\delta$-hyperbolic space, then we can choose $\Phi=\Phi_\delta$ and recall that there are $\Phi$-chains along geodesics in $M$. Since $\Phi_\delta$ is
determined from $\delta$, the $\Phi$-dependence of the estimates now reduces to a $\delta$-dependence.

\begin{corollary}[Green's function along hyperbolic geodesics]
\label{thm:gfalg} Assuming \refbac{} for $(M,L)$, $M$ $\delta$-hyperbolic, let $x,y,z\in M$ such that $y$ lies on geodesic connecting $x$ and $z$ with $d(x,y),d(y,z)>22\delta$.
Then there is a constant $c(\sigma,\ell, k,  n,\tau,\delta) >1$ such that
\[
    c^{-1}G(x,y)G(y,z)\le G(x,z)\le c\,G(x,y)G(y,z)\ .
\]
\end{corollary}

\subsection{Boundary Harnack Inequality}
\label{sec:bhp}
Using \thmref{thm:bhid} as a blueprint, we want to formulate a boundary Harnack inequality near points on the Gromov boundary of a $\delta$-hyperbolic space. 

As a replacement for balls in the classical version of the boundary Harnack inequality, we need some characterization of neighborhoods of a point at infinity. This is made precise by the notion of $\Phi$-neighborhood bases.
\begin{definition}[$\Phi$-neighborhood basis]
We call two open subsets $V\supset W$ of the space $M$ the \textbf{$\bm{\Phi}$-neighborhoods with hub $h\in M$}, if $\overline{W}\subset V$, $B_{\Phi_0}(h) \subset V\setminus\overline{W}$ and any two points $p \in \del V$ and $q \in \del W$ can be joined by
a $\Phi$-chain that has $h$ as a track point. We call an infinite family of nonempty open $\N_i\subset M$, $i=1,2,3,\dots$, with $\bigcap_i\N_i=\emptyset$ a \textbf{$\bm{\Phi}$-neighborhood basis}, if $\N_i$ and $\N_{i+1}$ are $\Phi$-neighborhoods with hub $p_i$, for every $i$.
\end{definition}
Just as metric balls are, besides their rôle as neighborhood bases, the basic playground for Harnack inequalities, these $\Phi$-neighborhoods are their counterpart in \emph{boundary} Harnack inequalities.

In $\delta$-hyperbolic spaces, every point in the Gromov boundary has a canonical neighborhood basis that is also a $\Phi_\delta$-neighborhood basis. Namely, as shown in \cite[Proposition 8.10]{bhk01}, we have

\begin{lemma}[$\Phi_\delta$-neighborhood basis] \label{lemma:den} For a $\delta$-hyperbolic manifold $M$,  $z \in \delG M \subset\overline{M}^{\mathrm{G}}$ and $\gamma: \R \ra  M$ a geodesic so that $z$ is its endpoint for $t \ra -\infty$ we set
\[\N^\delta_i:=\{x \in X\,|\, \dist\left(x,\gamma([300 \cdot i,+\infty))\right)<\dist\left(x,\gamma((-\infty,300 \cdot i])\right)\}\]
	Then the $\N^\delta_i$ are a $\Phi_\delta$-neighborhood basis and their closures  $\overline{\N^\delta_i} \subset \overline{M}^{\mathrm{G}}$ a neighborhood basis of $z \in\delG M$.
\end{lemma}

For non-self-adjoint operators (which will be important for the applications in \autoref{sec:b}) it is not always possible to find $L$-harmonic functions that vanish at infinity, because even minimal Green's functions might diverge. Hence we introduce a more general notion which can be thought of as minimal growth condition. This will be further explained in \thmref{prop:LvanM} below.

\begin{defprop}[\textit{L}-vanishing]\label{Lvan}
We say that a positive $L$-superharmonic function $u$ \emph{\textbf{$\bm{L}$-vanishes} towards infinity} on a domain $V\subset M$, if one of the following equivalent conditions is satisfied:
\begin{enumerate}
\item There is a positive $L$-superharmonic function $w$, such that $u/w\to0$ at infinity, i.e., for every $\epsilon>0$ there is a compact set $K\subset M$ with $u/w<\epsilon$ on $V\setminus K$.
\item There is a potential $p$ such that $p\ge u$ on $V$.
\item  The reduit $\red_u^V$ is a potential.\footnote{\cite{anc87} uses the first definition, but it is easier to employ the last.}
\end{enumerate}
\begin{proof}
(i)$\Rightarrow$(iii): Assume there is a positive $L$-harmonic function $h$, such that $\red_u^V\ge h>0$ on $M$. For some fixed $\epsilon>0$, choose a compact set $K$ with $u<\epsilon w$ on $V\setminus K$. By the properties of the reduit, we even have
\[\epsilon w\ge\red_u^{V\setminus K}\ge\red_u^V-\red_u^{V\cap K}\ge h-\red_u^{V\cap K}
\]
on all of $M$. Now $\red_u^{V\cap K}$ is a potential since $V\cap K$ is relatively compact in $M$, $\epsilon w-h$ is $L$-superharmonic and $\epsilon w\ge u=\red_u^V\ge h$ on $\del(V\setminus K)=\del(M\setminus(V\setminus K))\subset M$, thus we can apply the \thmnameref{gmp} to see $\epsilon w\ge h$ on all of $M$. Since $\epsilon$ was arbitrary, $h=0$.
\medskip

(iii)$\Rightarrow$(ii): Choose $p=\red_u^V$.
\medskip

(ii)$\Rightarrow$(i): It suffices to show that a potential $p$ satisfies the condition everywhere. Towards this end, consider the functions $\red_p^{M\setminus \overline{B_R(x_0)}}$ on balls around an arbitrary basepoint $x_0\in M$. They converge to zero for $R\to\infty$ because the limit is $L$-harmonic and $\le p$. Let $(x_j)$ be a countable dense set in $M$. We may choose a sequence $(R_i)$ such that $\red_p^{M\setminus \overline{B_{R_i}(x_0)}}(x_j)\le 2^{-i}$ for all $j\le i$. Then the function $w=\sum_i\red_p^{M\setminus \overline{B_{R_i}(x_0)}}$ is finite on a dense set, $L$-superharmonic and $p/w\le 1/i$ outside of $B_{R_i}(x_0)$.
\end{proof}
\end{defprop}

Note that all potentials such as the minimal Green's function are $L$-vanishing on $V=M$ and hence on all open sets because the property of $L$-vanishing is conserved on subsets as can be easily seen using condition (i).

On bounded sets, $L$-vanishing \emph{at infinity} is trivially true for any $L$-superharmonic function. The concept is also useless for unbounded sets shrinking to a small set, e.g.\ a single point,  when approaching infinity. It becomes significant for $V=W \cap M$ where $W\subset \overline{M}^\mathrm{G}$ is open with non-empty intersection $W \cap \delG M$. For more elaborate criteria in the context of Martin theory see \thmref{prop:LvanM}.
\medskip

On $\Phi$-neighborhoods we can now formulate the following central result.
\begin{theorem}[boundary Harnack inequality]
\label{thm:bhi} Assume $(M,L)$ satisfies assumptions \refbac.
Let $V\supset W$ be $\Phi$-neighborhoods with hub $h$ and $u$, $v$ two positive $L$-superharmonic functions that are $L$-harmonic and $L$-vanishing on $V$, then there is a constant $H_B=H_B(\sigma,\ell, k, n,\tau,\Phi)$ such that
\[
    \frac{u(x)}{u(y)}\le H_B \frac{v(x)}{v(y)}\quad\text{for any $x,y\in W$}\ .
\]
\begin{proof}
By \thmref{Lvan}, the reduit $\red_u^{V}$ is a potential and therefore admits a representation as
\[ \red_u^{V}=\int_{\del V} G(\cdot,z)\d\nu(z) \]
for some Radon measure $\nu$. On $\overline{W}$, $\red_u^{V}$ agrees with $u$ and therefore we have
\[
  u(x)=\int_{\del V} G(x,z)\d\nu(z)\le c\int_{\del V} G(x,h)G(h,z)\d\nu(z)
  = c\,G(x,h)u(h)\quad\text{for $x\in\del W$}
\]
using the assumption that $x\in\del W$ and $z\in\del V$ can be connected by a $\Phi$-chain through $h$ and \thmref{thm:gphi}. The other inequality from there gives
\[
  v(x)\ge c^{-1}\,G(x,h)v(h)\quad\text{for $x\in\del W$}
\]
and both inequalities extend to $W$ by the \thmnameref{gmp} because $\red_u^{V}$ and $G(\cdot,h)$ respectively are potentials. We can combine them to obtain
\[ \frac{u(x)}{v(x)}\le c^2\frac{u(h)}{v(h)}\quad\text{for $x\in W$}\ .\]
Interchanging the rôles of $u$ and $v$ yields the result with $H_B=c^4$.
\end{proof}
\end{theorem}

In most applications, $u$ and $v$ are either globally $L$-harmonic functions or minimal Green's functions with pole outside of $V$.
\medskip

In the case of a $\delta$-hyperbolic manifold the size of the smaller neighborhood can be explicitly quantified using \thmref{lemma:den}. We get

\begin{corollary}[hyperbolic boundary Harnack inequality]
\label{cor:hbhi} If $(M,L)$ satisfies \refbac{} and $M$ is $\delta$-hyperbolic, there is some positive constant $H_B(\sigma,\ell, k, n,\tau,\delta)>1$ such that two positive $L$-superharmonic functions $u$, $v$ that are $L$-harmonic and $L$-vanishing on a $\Phi_\delta$-neighborhood $\N^\delta_i$ of $\xi\in\delG M$ satisfy
\[ \frac{u(x)}{u(y)}\le H_B \frac{v(x)}{v(y)}\quad\text{for any }x,y\in \N^\delta_{i+1}\ . \]
\end{corollary}

\subsection{The Martin Boundary}\label{sec:mb}
In the previous section we have already formulated \emph{boundary} Harnack inequalities on the complete manifold $M$ although they neither needed not referred to any concrete boundary.
The natural boundary concept related to these boundary inequalities we introduce here is the  Martin boundary. These ideal boundaries can equally be defined on non-complete manifolds like Euclidean domains.

The Martin boundary is determined from both the geometry of $M$ and the analysis of $L$. A basic problem in Martin theory is the characterization of this boundary with a minimal amount of analytic input.  The distinguished result in our case is that when we also assume Gromov hyperbolicity, then the Martin boundary actually only depends on the geometry. It is homeomorphic to the Gromov boundary $\delG M$ of $M$.

We first recall some basic notions from Martin theory, cf.~\cite{anc90} or \cite[7.1]{pin95}.

\begin{description}
\item[Martin Boundary] \label{mb} For a non-compact Riemannian manifold $M$ and a linear
 second order elliptic operator $L$ on $M$ with a minimal Green's function $G: M\times M \ra (0,\infty]$ as well as a basepoint $p \in M$, we consider the space $S$ of
 sequences $s$
 of points $p_i \in M$, $i=1,2,\dots$
 with \begin{itemize}
    \item $p_i \in M$ has no accumulation points in  $M$,
    \item $K_{p_i}:= G(\cdot,p_i)/G(p,p_i) \overset{n\to\infty}{\longrightarrow} K_s$ compactly, for some function $K_s$.
 \end{itemize}
 The  \textbf{Martin boundary} $\delM(M,L)$  is the quotient of $S$ modulo the relation $s \sim s^*$ if $K_s \equiv K_{s^*}$. For $\zeta\in \delM(M,L)$, this function is written $K_\zeta$ and called the \textbf{Martin kernel}.
\end{description}

The Martin boundary does not depend on the choice of the basepoint $p$.  The Harnack inequality and elliptic theory show that each $K_\zeta \in \delM  (M,L)$ is a positive solution of $Lu = 0$ on $M$. This also shows
that the convex set $S_L(M)$ of positive solutions of $L u =0$ on $M$ with $u(p) =1$ is compact in the topology of compact convergence. In turn, $\delM  (M,L)$ is a compact
subset of $S_L(M)$.

The \textbf{Martin topology} on  $\overline{M}^\mathrm{M}:= M \cup \delM(M,L)$ is the topology of compact convergence on the space of associated Martin functions $\{K_y\:|\:y\in\overline{M}^\mathrm{M}\}$.

Then the usual topology is induced on $M\subset\overline{M}^\mathrm{M}$, $\delM(M,L)$ is closed and $\overline{M}^\mathrm{M}$ is compact. The space $\overline{M}^\mathrm{M}$ is called the \textbf{Martin compactification}. It is easy to see that $\overline{M}^\mathrm{M}$ is metrizable, cf.~\cite[Ch.\,I.7]{bj06} or \cite[Ch.\,12]{hm14} for further details.

To motivate the idea of Martin integrals we recall the following classical result, cf.~\cite[Ch.\,6]{cho69}:

\begin{proposition}[Minkowski's theorem]\label{mit} Each point in a convex set  $K \subset \R^n$ is a
convex combination of the extremal points of $K$.
\end{proposition}

\begin{figure}[h]
\centering
\includegraphics[width=0.84\textwidth]{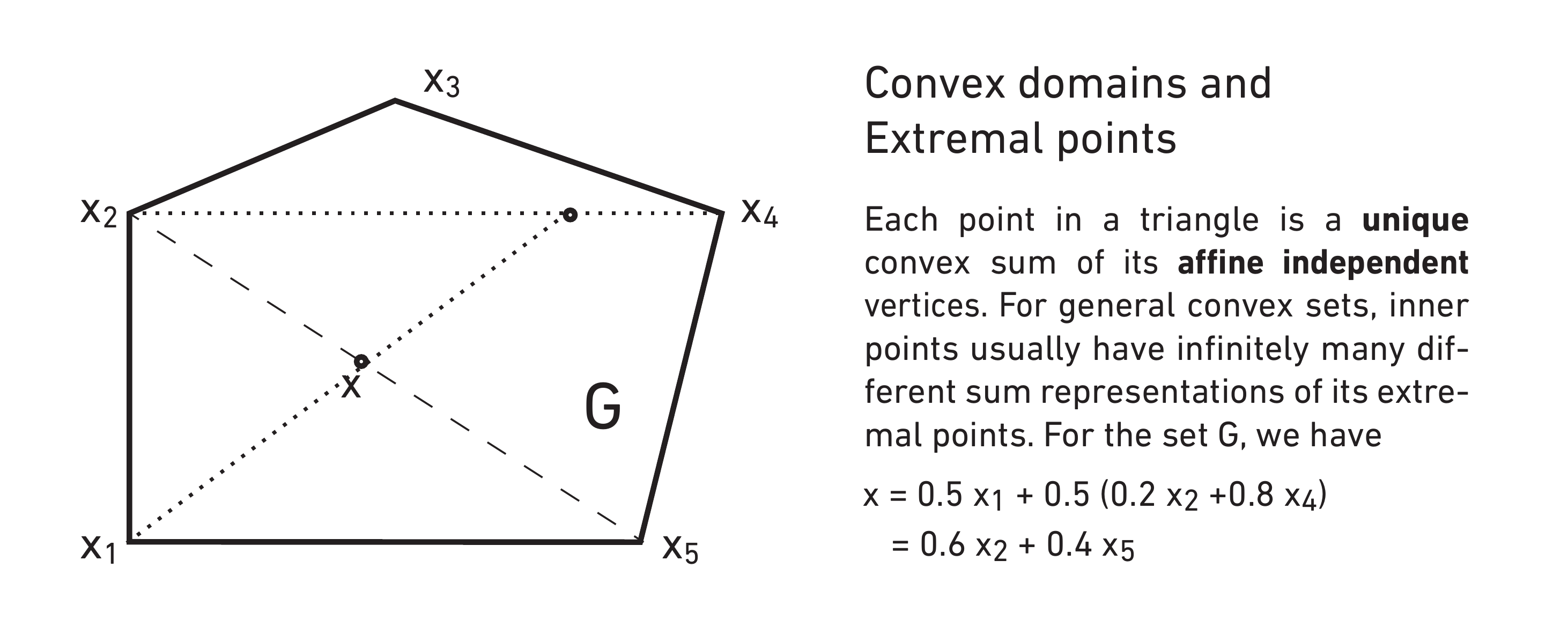}
\caption{Minkowski's Theorem in $\R^2$}
\end{figure}

The \textbf{Martin integral} is essentially an extension of this result to the case of the convex set $S_L(M)$. The extremal elements of $S_L(M)$ form a subset $\delM^0 (M,L) \subset \delM(M,L)$ of the Martin boundary one can think of as the vertices of an infinite dimensional simplex spanning $S_L(M)$. A positive solution $u$ of $Lu = 0$ on $M$ with $u(p)=1$ is extremal if and only if $u$ is a minimal solution. Here we call $u$ \textbf{minimal}
if for any other solution $v >0$,  $v \le u$, we have $v \equiv c \cdot u$, for some constant $c
>0$. Therefore $\delM^0 (M,L) \subset \delM (M,L)$ is also called the \textbf{minimal Martin boundary}.

The Choquet integral representations in \cite[Ch.\,6]{cho69} give the following general version of the Martin representation theorem, cf.~\cite[7.1]{pin95}:

\begin{proposition}[Martin integral]\label{thm:mrt} For any positive solution $u$ of $Lu = 0$ on $M$, there is a
unique Radon measure $\mu_{u}$ on $\delM^0 (M,L)$ so that
\[u(x)  =\int_{\delM^0 (M,L)} K_\zeta(x) \, \d \mu_u(\zeta)\ .\]
Conversely, for any Radon measure $\mu$ on $\delM^0 (M,L)$,
\[u_\mu(x)  =\int_{\delM^0 (M,L)} K_\zeta(x) \, \d \mu(\zeta)\]
defines a positive solution of $Lu_\mu = 0$ on $M$.
\end{proposition}

Although this already looks like a classical contour integral the result is not yet truly satisfactory. Unlike the classical case the boundary $\delM^0 (M,L)$ depends not only on the underlying space but also on the analysis of the operator $L$. A natural question is whether one could get rid of this dependence. In general the answer is no. In \ref{ex:nun} we will see some by no means exotic examples. However, we will now see that in our case of operators on $\delta$-hyperbolic spaces, this is actually possible. This is a remarkable particularity not even valid for such simple spaces as in \thmref{ex:pro}.

\paragraph{Back to hyperbolic spaces} In the situation at hand, $\Phi$-neighborhood bases are essentially neighborhood bases of minimal Martin boundary points.
\begin{theorem}[characterization of minimal Martin points]\label{cmmp}
Assume \refbac{} holds.
Let $(\N_i)$ be a $\Phi$-neighborhood basis with hubs $p_i$.
Denoting the interior of the closure of $\N_i\subset M \subset \overline{M}^\mathrm{M}$ in the Martin compactification $\overline{M}^\mathrm{M}$ of $M$ by $\tilde{\N}_i$, there is exactly one Martin boundary point $\zeta$ in $\bigcap\tilde{\N}_i$. The resulting $K_\zeta$ is characterized as the only positive $L$-harmonic function $L$-vanishing on every $M\setminus\N_i$ except for scalar multiples. In particular, this Martin point is minimal.

\begin{proof}
By the Harnack inequalities, the sequence $K_{p_i}=\frac{G(\cdot, p_i)}{G(p,p_i)}$ has a subsequence compactly converging to some $L$-harmonic function $K_\zeta$ representing a Martin boundary point $\zeta$.
$K_\zeta$ is $L$-vanishing on every $M\setminus\N_i$ because by the \thmnameref{thm:bhi}, every $K_{p_j}$ for $j\ge i$, and hence the limit, is upper bounded by the potential $H_B K_{p_i}$ on $M\setminus\N_i$.

Now assume there is another positive $L$-harmonic function $u$ that is $L$-vanishing on $M\setminus\overline{\N}_i$ for every $i$, w.l.o.g.\ $u(p)=1$. Applying the \thmnameref{thm:bhi} we see $H_B^{-1}K_\zeta\le u \le H_B K_\zeta$ on $M$.
Thus $\eta:=\inf u/K_\zeta\le1$ is positive. By the strong maximum principle \cite[Theorem 3.2.6, p.\,84]{pin95}, the $L$-harmonic function $u-\eta K_\zeta\ge0$ has to be positive everywhere, else it would be identical zero.
In the former case we can again apply the \thmnameref{thm:bhi} to $K_\zeta$ and $u-\eta K_\zeta$ to get $(\eta+(1-\eta)H_B^{-1})K_\zeta\le u$ which contradicts the definition of $\eta$, unless $\eta=1$ and $u=K_\zeta$.
\end{proof}
\end{theorem}

We have the following characterization of the Martin boundary in case we have enough $\Phi$-chains:
\begin{corollary}[identifying the Martin boundary]\label{imb}
If under assumptions \refbac{} in a given compactification $\overline M$ of $M$ (i.e., $\overline M$ is compact and $M\subset\overline M$ dense) every boundary point admits a neighborhood basis of the form $\overline\N_i\subset \overline M$ for some $\Phi$-neighborhood basis $(\N_i)$, it is canonically homeomorphic to the Martin compactification $\overline M^\mathrm{M}$.
\begin{proof}
\thmref{cmmp} yields an injective map from $\overline M$ to $\overline M^\mathrm{M}$. It is continuous because for every sequence $(y_i)$ in $\overline M$ converging to $\zeta\in\overline M\setminus M$ the corresponding Martin functions $K_{y_i}$ converge to the unique Martin function that is $L$-vanishing on all $M\setminus\N_i$ for some $\Phi$-neighborhood basis $(\N_i)$ of $\zeta$, that is $K_\zeta$. Thus, by elementary properties of compactifications, it is already a homeomorphism.
\end{proof}
\end{corollary}

From this and \thmref{lemma:den} we get the following principal potential theoretic result on Gromov hyperbolic manifolds.
\begin{corollary}[Gromov boundary and Martin boundary] \label{cor:gmr}   Assume that $M$ is Gromov hyperbolic and \refbac{} holds. Then the Gromov and Martin boundaries of $M$ are canonically homeomorphic and every Martin boundary point is already minimal,
\[ \delG M\cong \delM (M,L)\cong  \delM^0 (M,L) .\]
In particular, some function $u >0$ on $M$  solves $Lu = 0$ if and only if there is a (unique) Radon measure $\mu_u$ on $\delG M$
such that
\[u(x)  =\int_{\delG M} K_\zeta(x) \, \d \mu_u(\zeta).\]
\end{corollary}
The validity of such a simple identification of the Martin boundary with a geometric boundary, which one may possibly expect from a naive guess, actually is a rare exception. If we only slightly violate the hyperbolicity constraint we can get a completely different and rather inscrutable outcome with many non-minimal Martin boundary points.
\begin{example}[Ideal Boundaries of $\mathbb{H}^m \times \mathbb{H}^n$] \label{ex:pro}
The product space of two classical hyperbolic spaces  $\mathbb{H}^m \times \mathbb{H}^n$, $m, n \ge 2$, has bounded geometry but it is no longer Gromov hyperbolic since it contains flat planes obtained as products of pairs of geodesics in the two factors. In \cite{gw93} we find a thorough discussion of the case of the Laplace operator.  There exist positive functions $h$ with $-\Delta h = \lambda \cdot h$ if and only if $\lambda \le \lambda_0$ where $\lambda_0$ is the generalized principal eigenvalue. For $\lambda<\lambda_0$, the operator $-\Delta - \lambda \id$ is adapted and weakly coercive.
The boundary at infinity (a natural generalization of Gromov boundary) is homeomorphic to $S^{n+m-1}$, cf.\ \cite[II.8.11(6), p.\,266]{bh99}. In turn, for the minimal Martin boundary of $-\Delta - \lambda\id$, for $\lambda<\lambda_0$,  we have
\[\delM^0 (\mathbb{H}^m \times \mathbb{H}^n,-\Delta - \lambda\id) = S^{m-1} \times  S^{n-1} \times I_\lambda\ ,\]
where $I_\lambda$ is a closed interval with a natural parameterization which depends on $\lambda$ and which degenerates to a single point when $\lambda \ra \lambda_0$ \cite[p.\,21]{gw93}. The full Martin boundary contains two additional pieces
\begin{equation}\label{part}
\delM (\mathbb{H}^m \times \mathbb{H}^n,-\Delta - \lambda\id) = (S^{m-1} \times \mathbb{H}^n) \cup (S^{m-1} \times  S^{n-1} \times I_\lambda) \cup (\mathbb{H}^m \times S^{n-1})/\sim\ .
\end{equation}
The gluing maps for $\sim$ are described in \cite[p.\,27]{gw93}. We observe that not only the boundary at infinity does not coincide with $\delM^0 (\mathbb{H}^m \times \mathbb{H}^n,-\Delta - \lambda\id)$ but the details of the partition of the full Martin boundary (\ref{part}) depend on $\lambda$.\qed
\end{example}

Now that we have the right notion for a potential theoretic boundary at infinity, there are some new formulations for $L$-vanishing. Note that these are slightly different from the formulation in \thmref{Lvan} because there it was only possible to refer to $L$-vanishing \emph{at infinity of open sets $V$ in $M$}, i.e., on the Martin boundary points in $\overline{V}\cap\delM M\subset\overline{M}^\mathrm{M}$.
\begin{proposition}[\textit{L}-vanishing and Martin boundary]\label{prop:LvanM}
Assume \refbac{} holds and every Martin boundary point has a $\Phi$-neighborhood basis, e.g., $M$ is Gromov hyperbolic.
For an open subset $\varXi\subset \delM M$ of the Martin boundary and a positive $L$-harmonic function $u$ on $M$ the following are equivalent:
\begin{enumerate}
\item $u$ is $L$-vanishing on any open set $V\subset M$ with $\overline V\cap\delM M\subset\varXi$ in the Martin compactification.
\item On any open set $V\subset M$ with $\overline V\cap\delM M\subset\varXi$ in the Martin compactification, the following property holds: Each positive $L$-harmonic function $v$ on $V$ with $v\ge u$ on $\del V\cap M$ satisfies $v\ge u$ on $V$.
\item The Martin measure $\mu_u$ associated to $u$ is supported outside $\varXi$, i.e., $\mu_u(\varXi)=0$.
\end{enumerate}
\begin{proof}
(i)$\Rightarrow$(ii)
If $u$ is $L$-vanishing on $V$, $\red_u^V$ is a potential and $\red_u^V=u$ on $V$. Then the \thmnameref{gmp} gives exactly the desired property.
\medskip

(iii)$\Rightarrow$(i) Each Martin function $K_{\zeta}$ is $L$-vanishing outside $\zeta$ and, therefore, the Martin integral representing $u$ $L$-vanishes outside the support of  $\mu_u$.
\medskip

(ii)$\Rightarrow$(iii)
Assume $\mu_u(\varXi)\neq0$, then there is a compact $K\subset\delM M$ and an open $W\subset\overline{M}^\mathrm{M}$ such that $K\subset W\cap\delM M\Subset\varXi$, $V:=W\cap M$ satisfies the condition in (ii), and $\mu_u(K)\neq0$ (since the Radon measure $\mu_u$ is inner regular). Therefore it is enough to consider the case where $u \equiv u_K:=\int_K K_\zeta\,\d\mu_u(\zeta)$.

We compare $u_K$ with the minimal Green's function $G(\cdot,p)$ with pole $p\in M\setminus V$. Recalling the argument of (iii)$\Rightarrow$(i) we know that $u_K$ is $L$-vanishing on an open neighborhood $N$ of $M\setminus V$. By compactness of $\overline{M\setminus V}$ in the Martin compactification, $M\setminus V$ can be covered by finitely many $\Phi$-neighborhoods contained in $N$ and a compact subset of $M$. Then the \thmnameref{thm:bhi} shows that there is a $C>0$ such that $C\cdot G(\cdot,p)\ge u_K$ on $M\setminus V$ and especially on $\del V$. But then from (ii) it follows that $C\cdot G(\cdot,p)\ge u_K$ on $V$, hence on all of $M$, hence $u_K\equiv 0$ because $G(\cdot,p)$ is a potential.
\end{proof}
\end{proposition}

In the next section we will see that the integral representation we have obtained in the Gromov hyperbolic case naturally extends classical contour integrals for (positive) harmonic functions on the subclass of Euclidean domains which are \emph{uniform}.

\section{Geometry and Analysis on Uniform Domains}\label{sec:b}

Here we consider domains sharing basic geometric and analytic properties with smoothly bounded ones. Examples are canonical deformations  of such domains to Gromov hyperbolic manifolds of bounded geometry but also the validity of Poincar\'{e}-Sobolev type inequalities and embedding theorems. We start with the definition and examples of so-called uniform domains. Then we explain how one can transfer Ancona's theory from hyperbolized uniform domains back to original Euclidean domains to derive new analytic results and to recover classical integral formulas.

\subsection{Uniform Domains}

The path connectedness of an open set is a topological condition. For finer geometric and analytic investigations one seeks for a quantitative form of path connectedness. One of the nowadays central notions is described in the following definition. We refer to \cite{aik12} for an instructive overview and comparison of other regularity concepts for domains.

\begin{description}
\item[Uniform Domains]\emph{A domain $D\subset\R^n$ is called a \textbf{uniform domain}, more precisely, a $c$-uniform domain, if there is a $c\ge 1$ such that any two points $p,q \in D$ can be joined by a
\textbf{c-uniform curve}. That is a rectifiable path $\gamma: [a,b] \to D$, for some $a<b$, from $p$ to $q$ so that the following conditions are satisfied:
\begin{itemize}
  \item \emph{\textbf{Quasi-geodesic:}} $l(\gamma)\le c\, d(p,q)$.
  \item \emph{\textbf{Twisted double cones:}} $\min\{l(\gamma|_{[a,t]}),l(\gamma|_{[t,b]})\} \le c \, \dist(\gamma(t),\del D)$ for any $t\in[a,b]$.
\end{itemize}}
\end{description}

Note that being $c$-uniform is a scaling invariant condition: whenever $D$ is $c$-uniform,  $\lambda \cdot D$, $\lambda >0$, is also $c$-uniform.

\begin{examples} [uniform and non-uniform domains] \label{ex:nud}  We start with some types of Euclidean domains which are uniform:
\begin{itemize}
\item Any bounded domain with smooth or at least Lipschitz boundary is uniform.
\item Non-compact rotationally symmetric domains bounded with profile functions of at least linear growth are uniform. As an explicit example, choose $f(t)= c_1 \cdot t + c_2,$ for constants $c_i >0$ and let $F:\R^{n-1}\ra \R$, $n\ge 3$, be given by $F(x):=f(|x|)$. Consider the domain \[D_f:=\{(x_1,..,x_n) \in \R^n \,|\, |x_1|< F(x_2,..,x_n) \}.\]
   To explain the idea of how to define the quasi-geodesics and  twisted double cones joining, e.g., the pairs of points $p^\pm_k=(0,...0,\pm k)$, $k \ge 1$, we note that for $k=1$ we can choose any such twisted double cone $\subset D_f$ along a half circle in the $n-1,n$-plane joining  $p^-_1$ and $p^+_1$. Then the quasi-geodesic and twisted double cone scaled by $k$ serve for the $p^-_k$ and $p^+_k$. Here it essential that $f$ grows at least linearly to ensure that the twisted double cone remains in $D_f$.
\item Bounded domains with certain types of fractal boundaries, like the Snowflake in $\R^2$ or the complement of the Sierpinski gasket in $\R^3$ \cite{aik03}, are uniform domains.
\end{itemize}
\begin{figure}[h]
\centering
\includegraphics[width=\textwidth]{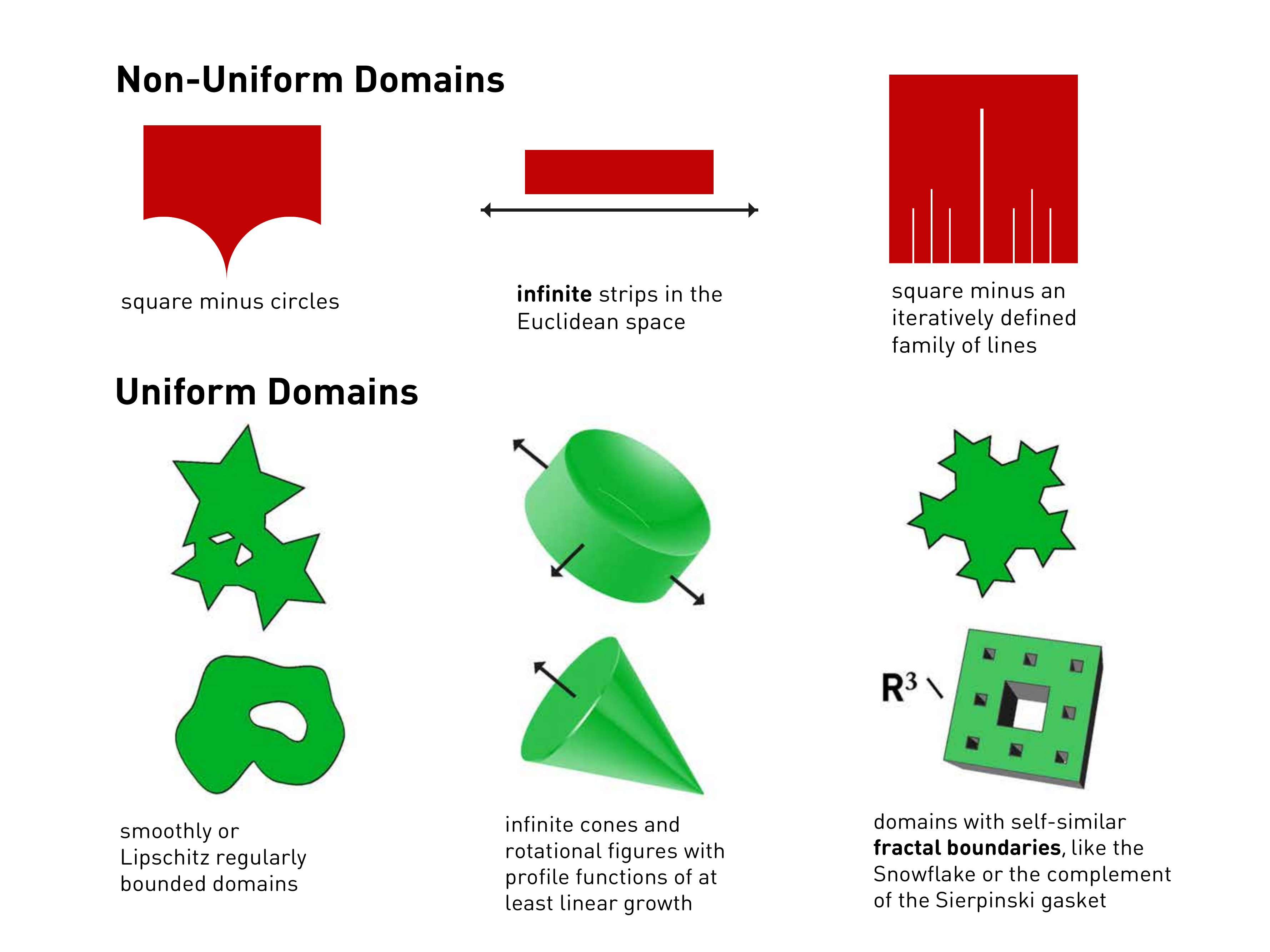}
\caption{Typical (Non-)Uniform Domains}
\end{figure}
However, even completely harmless domains can be non-uniform:
\begin{itemize}
\item The difference of the cube $(-1,1)^n \subset \R^n$ minus the ball $B^n_1(0) \subset \R^n$: $(-1,1)^n \setminus B_1(0)$ is not a uniform domain, since we cannot reach points arbitrarily near to the boundary point $(0,\dots,0,1)  \in \p \left((-1,1)^n \setminus B_1(0)\right)$ by twisted cones, for some common constant $c>0$.
\item The cylinder $B^{n-1}_1(0) \times \R \subset \R^n$  and similarly $B^k_1(0) \times \R^{n-k} \subset \R^n$, for $1 \le k \le {n-1}$ are non-uniform. For the pairs of points $p^\pm_k=(0,...0,\pm k)$, $k \ge 1$, there are only twisted double cones reaching from $p_k^+$ to $p_k^-$, for $c =c(k) \le 1/k$.
\item Similarly, the domain $\R^3 \setminus \R \times \Z^2$  is not uniform.
\item  As in the earlier example  of non-compact rotationally symmetric domains $D_f$ we choose profile functions but this time of sublinear growth like $f(t)= c_1 \cdot \sqrt{t} + c_2,$ for constants $c_i >0$. Then, as in the last counterexample,  $D_f$ is non-uniform.
\end{itemize}
\end{examples}

\subsection{Quasi-Hyperbolic Geometry}
On any domain $D\subset\R^n$ we can define the \emph{quasi-hyperbolic metric} $k_D$ introduced as an extension of the definition of the Poincar\'{e} metric on the unit disc to arbitrary domains. It is given by
\[
    k_D(x,y):=\inf\left\{\int_\gamma1/\dd\:\Bigg|\:\text{$\gamma\subset D$ rectifiable curve joining $x$ and $y$}\right\}
\]
where we set $\dd=\dist(\cdot,\del D)$. This corresponds to a conformal deformation of $g_\eucl$ to the merely Lipschitz continuous $\dd^{-2}g_\eucl$. The quasi-hyperbolic metric $k_D$ always has bounded geometry since the chart
\[
 \id:(B_{\dd(x)/2}(x),k_D)\to(B_{\dd(x)/2}(x),\dd^{-2}(x)g_\eucl)
\]
   (with radii measured in Euclidean distance) is 2-bi-Lipschitz and we see that $(D,k_D)$ has  $(\sigma,\ell)$-bounded geometry for some $\sigma,\ell>0$ both independent of $D$. \\

In general, $k_D$ need not to be (Gromov) hyperbolic, except $D$ is a uniform domain. Due to geometric work, in particular of Bonk, Heinonen, Koskela \cite{bhk01}, we have:

\begin{theorem}[Gromov uniformization]\cite[Th.\,1.11]{bhk01} \label{thm:grh} A bounded  domain   $D \subset \R^n$  is  uniform  if  and only  if  $(D,k_D)$ is both  Gromov  hyperbolic and the
Euclidean boundary $\p D$ of $D$ is naturally  quasisymmetrically  equivalent to  the  Gromov  boundary $\delG(D,k_D)$. Moreover, when $D$ is $c$-uniform, then $(D,k_D)$ is $\delta_c$-hyperbolic of $(\sigma,\ell)$-bounded geometry, for some $\delta_c>0$ depending only on $c$.
\end{theorem}

A map $f: X \to Y$ between metric spaces $(X,d_X)$ and $(Y,d_Y)$ is called  \emph{quasi-symmetric} if $f$ is not constant and if there
is a homeomorphism $\eta: [0,\infty) \to [0,\infty)$, so that for any $x, a, b\in X$, $t\ge 0$, and $d_X(x,a)\le  t\cdot d_X(x,b)$ it follows that $d_Y(f(x),f(a))\le  \eta(t)\cdot d_Y(f(x),f(b))$.
\medskip

Since $\dist(\cdot,\del D)$  is Lipschitz continuous but not smooth, $(D,k_D)$ is not a Riemannian manifold. But we need some smoothness to ensure the adaptedness of $L$, in particular to get uniform bounds on the coefficients of the adjoint operator $L^*$. Thus, in particular for analytic considerations, we will work with smoothed versions of $\dd$. Fortunately, there is a well-controlled smoothing of $\dd$:
\begin{theorem}[Stein-Whitney smoothing]\cite[VI.2.1, p.\,171]{ste70} \label{thm:stw} There are constants $c_\alpha=c_\alpha(n)$ for every multiindex $\alpha$ such that on every open set $D\subset\R^n$ there is a $C^\infty$ function $\tdd$ defined on $D$ with
\[c_0^{-1}\dd\le\tdd\le c_0\dd\quad\text{and}\quad \left|\diffp[\alpha]{\tdd}{x}\right|\le c_\alpha \dd^{1-|\alpha|}\ .\]
In particular, $(D,\dd^{-2}g_\eucl)$ and $(D,\tdd^{-2}g_\eucl)$ are quasi-isometric and, hence,  $(D,\dd^{-2}g_\eucl)$ is Gromov hyperbolic if and only if $(D,\tdd^{-2}g_\eucl)$ is Gromov hyperbolic. Moreover, when $(D,\dd^{-2}g_\eucl)$ is $\delta$-hyperbolic, for some $\delta>0$, then $(D,\tdd^{-2}g_\eucl)$ is $a_n \cdot \delta$-hyperbolic, for some $a_n \ge 1$ depending only on $n$.
\end{theorem}

\subsection{Martin Theory for (Non-)Uniform Domains}
\label{sec:nu}

Quasi-hyperbolic metrics can be used to transfer  the potential theory on Gromov hyperbolic manifolds to understand the corresponding theory on uniform domains.
Since there are examples of non-uniform domains which carry a rather complicated potential theory, this correspondence also gives us non-hyperbolic manifolds of bounded geometry with a much less transparent Martin theory than in the hyperbolic case.

On a uniform domain $D \subset \R^n$, we consider  a Schrödinger operator $L^\eucl=-\Delta+V$ with some smooth potential $V \le a \cdot \tdd^{-2}$, for some $a>0$. That is, we even allow a controlled diverging behavior of $V$ towards $\p D$. The methods apply to more general operators, but we only want to explain the idea and this case is particularly important and notationally simple.

On  $(D,\tdd^{-2}g_\eucl)$ we consider the operator $L=\tdd^2 \cdot L^\eucl$. Then $L$ is adapted in the bounded geometry charts. Note that $L$-harmonic functions are the same as $L^\eucl$-harmonic functions, but $L$ is not self-adjoint anymore. At least the adjoint $L^*$ is still  adapted.

 $L$ is  weakly coercive relative to $(D,\tdd^{-2}g_\eucl)$ provided $L^\eucl$ has a positive first eigenvalue not only for  $L^\eucl$ but also
 for $\tdd^2 \cdot L^\eucl$ in the original Euclidean metric \cite{anc86}:
\begin{description}
\item[Strong Barrier] \emph{There is a function $s>0$ and some $\ve>0$ such that $L^\eucl s \ge \ve \cdot \tdd^{-2} \cdot s$.}
\end{description}

\begin{example}[regularity versus strong barrier] \label{ex:stba}  The Laplace operator on bounded domains with smooth or Lipschitz  boundary admits a strong barrier \cite{anc86,kad66}.
In general, the strong barrier condition is not a consequence of the uniformity of $D$ but of  additional exterior conditions. An exterior twisted cone condition is sufficient \cite{aik12}.
\end{example}

Under these mild conditions on $L^\eucl$ the potential theory we established above readily transfers to uniform domains where we now use the sets $\N^\delta_i$ from \thmref{lemma:den}, defined relative $\tdd^{-2}g_\eucl$, as a replacement for the concentric balls of \thmref{thm:bhid} in the classical Euclidean setup.

\begin{corollary}[Martin theory on uniform domains]
\label{cor:mtud} Let $L^\eucl=-\Delta+V$ be an operator with $V \le a\cdot \tdd^{-2}$ satisfying  a strong barrier condition on some $c$-uniform domain $D$, for some $a,\ve, c>0$. Then we have for any two $L^\eucl$-harmonic functions $u$, $v >0$ on $\N^\delta_i \cap D$ both $L^\eucl$-vanishing\footnote{The definition of $L^\eucl$-vanishing towards $\N^\delta_i \cap \p D$ is the same as in \thmref{Lvan} and \thmref{prop:LvanM} where we merely had an ideal boundary of $M$ an hence called it $L$-vanishing towards infinity.} towards $\N^\delta_i \cap \p D$,
\[
u(x)/v(x) \le C \cdot  u(y)/v(y) \mm{ \emph{for any two points} } x, y \in \N^\delta_{i+1} \cap D\ ,
\]
for some $C(a,\ve, c,n)\ge 1$.  Thus, the topological and the Martin boundary are homeomorphic and every Martin boundary point is  minimal:\,
$\del D \cong \delM^0(D,L^\eucl) \cong \delM(D,L^\eucl)$.\\

\begin{proof} We first apply  \thmref{thm:grh} to the uniform domain $D$, then we transform $\dd^{-2}g_\eucl$  into the smoothed $\tdd^{-2}g_\eucl$  using the \thmnameref{thm:stw}. We get an $\alpha_n \cdot \delta_c$-hyperbolic manifold of $(\sigma,\ell)$-bounded geometry with
\[
\del D \cong \delG (D,\dd^{-2}g_\eucl) \cong \delG (D,\tdd^{-2}g_\eucl)\ .
\]
Now we can study $L$, which is weakly coercive and adapted, via the \thmnameref{cor:hbhi}  and  \thmref{cor:gmr}. We conclude that
\[\del D \cong  \delG (D,\tdd^{-2}g_\eucl) \cong \delM ((D,\tdd^{-2}g_\eucl),L)\cong  \delM^0 ((D,\tdd^{-2}g_\eucl),L)\ .
\]
The Green's functions relative $L^\eucl$ on  $(D,g_\eucl)$ and $L$ on  $(D,\tdd^{-2}g_\eucl)$ are related by $G(x,y)=\tdd^{n-2}(y)G^\eucl(x,y)$, while the Martin functions and the solutions, along with the notion of $L$-vanishing, are the same. This means that the boundary Harnack inequalities carry over to $D$ with $L^\eucl$  and that $\delM^0((D,\tdd^{-2}g_\eucl),L)\cong \delM^0(D,L^\eucl)$ and $\delM ((D,\tdd^{-2}g_\eucl),L)\cong \delM(D,L^\eucl)$.
\end{proof}
\end{corollary}

Aikawa \cite{aik01,aik04} has used a somewhat different approach to prove some remarkable refinements underlining the sharpness of these potential theoretic result.
\begin{itemize}
  \item In \cite{aik01} he derives boundary Harnack inequalities for the Laplacian on arbitrary uniform domains even \emph{without} imposing  a strong barrier condition.
  However, in his result the Harnack constant depends on the domain $D$, whereas in the previous result it only depends on the two parameters $\ve$ and $c$.
  \item In \cite{aik04} we even find that $D \subset \R^n$ is a uniform domain if and only if the Laplacian satisfies boundary Harnack inequalities. As in \thmref{ex:stba}
 one also needs to assume an exterior regularity of $D$. Again an exterior twisted cone condition is sufficient \cite{aik12}.
\end{itemize}

\begin{figure}[h]
\centering
\includegraphics[width=\textwidth]{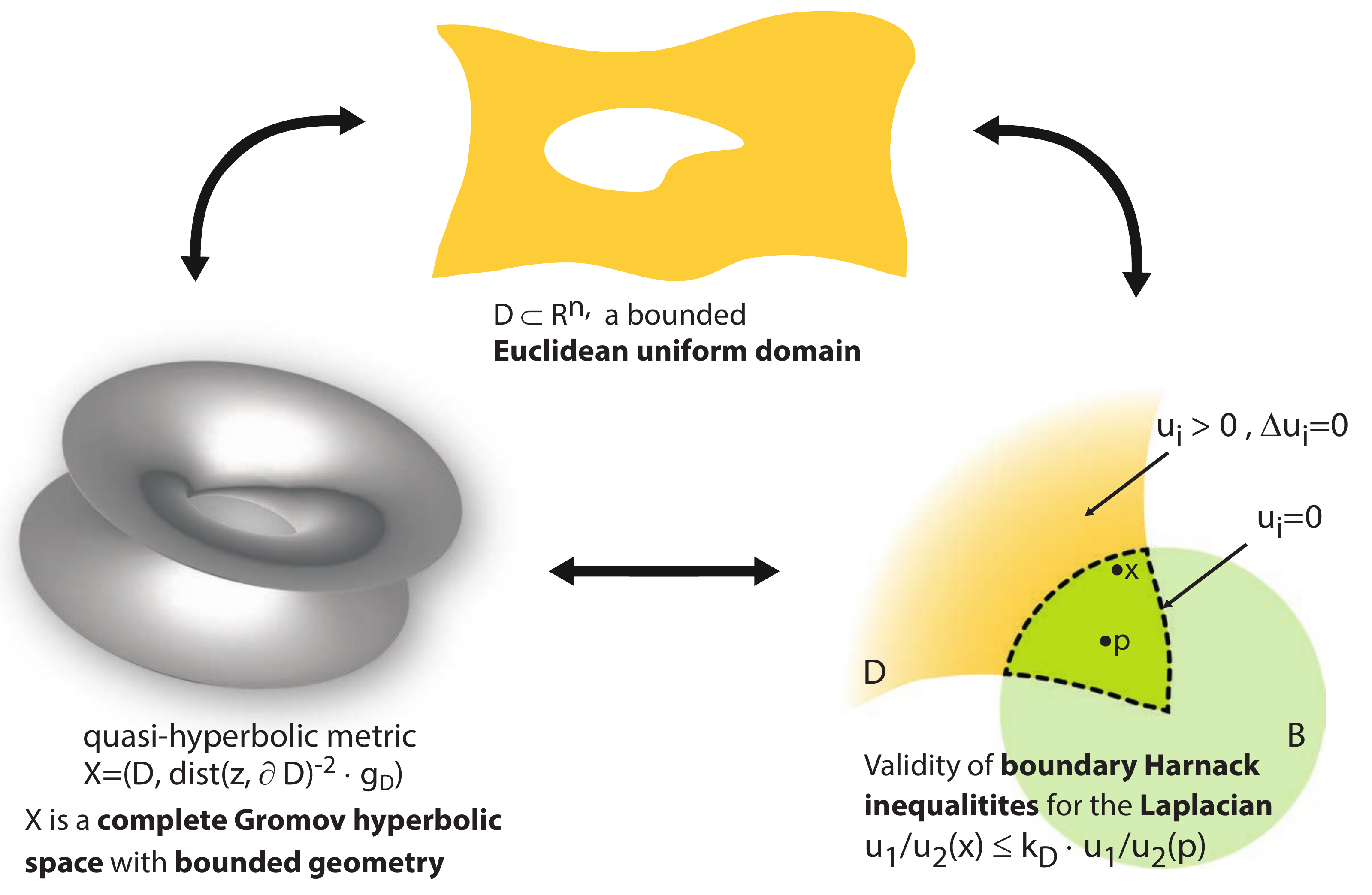}
\caption{The three essentially \textbf{equivalent} properties of a Euclidean domain $D$}
\end{figure}

\begin{examples} \label{ex:nun} We mention two instructive examples where non-uniformity of a domain destroys the existence of boundary Harnack inequalities and where we concretely see how far the topological boundary may
deviate from the Martin boundary of the Laplacian.
\begin{itemize}
  \item In  \cite{anc12} Ancona has described non-uniform Euclidean cones with only one topological point at infinity but with uncountably many minimal Martin boundary points at infinity.
  \item In \cite{IP94} Ioffe and Pinsky have proved that for the non-uniform rotationally symmetric domains $D_f \subset \R^n$ from \thmref{ex:nud} the
set of  Martin boundary points at infinity is homeomorphic to $S^{n-2}$.
\end{itemize}
\end{examples}

Finally we notice that this Martin theory on uniform domains reproves classical contour integral formulas, for instance, for harmonic functions. The Herglotz theorem \cite{her11,dur83} shows that a function $f >0$ on the Euclidean unit disk $(D,g_\eucl)$ is harmonic, $\Delta_\eucl \, f = 0$, if and only if there is a Radon measure $\mu_f$ on $S^1$ such that
\[
 f(x)=\int_{S^1} \frac{1-|x|^2}{|x-y|^2}\,\d\mu_f(y)\ .
\]
For the Green's function $G$ for $\Delta_\eucl$, and thus for the Green's function of the associated hyperbolized manifold and operator, a direct computation shows that \[K_\zeta(x) = \lim_{z \ra \zeta} \frac{G(x,z)}{G(0,z)}=\frac{1-|x|^2}{|x-\zeta|^2}, \mm{ for any } \zeta \in S^1.\]
But we also know from the above that we have a unique Martin integral representation:
\[
 f(x)=\int_{S^1} K_\zeta(x)\, \d\mu_f(\zeta)\ .
\]
 $\mu_f$ is now understood as a measure on a certain ideal boundary  of  $(D,g_\text{hyp})$,  the Martin boundary, which in this case equals $S^1$.

\phantomsection\addcontentsline{toc}{section}{\refname}
\bibliographystyle{aomalpha2}
\bibliography{lit2}

{\small
\textsc{Mathematisches Institut, Universit\"at M\"unster, Einsteinstra\ss e 62, Germany}

\textit{E-mail adresses:} \email{m.kemper@uni-muenster.de}, \email{j.lohkamp@uni-muenster.de}

\end{document}